\documentclass{article}
\usepackage[utf8]{inputenc}
\usepackage{inputenc}
\usepackage{amsmath}
\usepackage{amssymb}
\usepackage{ upgreek }
\usepackage{graphicx}
\usepackage{layout}
\usepackage{dsfont}
\usepackage{bbold}
\usepackage[english]{babel}
\usepackage[top=2cm, bottom=3cm, left=2cm, right=2cm]{geometry}
\usepackage{multicol}
\usepackage{amsthm}
\usepackage{mathrsfs}
\usepackage{tikz}
\usepackage{tikz-network}
\usepackage[toc,page]{appendix}
\usepackage{hyperref}
\usepackage{subfig}

\newtheorem{theorem}{Theorem}
\newtheorem{corollaire}[theorem]{Corollary}

\newtheorem{definition}[theorem]{Definition}
\newtheorem{conjecture}[theorem]{Conjecture}

\newtheorem{proposition}[theorem]{Proposition}

\newtheorem{lemma}[theorem]{Lemma}
\newtheorem{exemple}[theorem]{Example}
\newtheorem{remark}[theorem]{Remark}

\newcommand{\R}{\mathbb{R}} 
\newcommand{\N}{\mathbb{N}} 
\newcommand{\Z}{\mathbb{Z}} 
 
\newcommand{\al}{\alpha} 
\newcommand{\be}{\beta}

\newcommand{\si}{\sigma}

\newcommand{\bs}{\backslash} 
 
\newcommand{\eps}{\upvarepsilon}

\newcommand{\bip}{{\sc Bipartite Influence}}
\newcommand{\poscnf}{{\sc Pos Cnf}}

\newcommand{\rmv}[2]{
\mbox{Rmv}(#1, #2)
}

\newcommand{\decisionpb}[4]{
        \begin{minipage}{#4\textwidth}
                #1\\
                \emph{Instance:} #2\\ 
                \emph{Question:} #3
        \end{minipage}
}

\usepackage{comment}

\title{Bipartite instances of INFLUENCE}
\author{Eric Duchêne, Nacim Oijid and Aline Parreau }
\date{October 2020}

\tikzstyle{vertexW}=[circle,inner sep=0, minimum size =6 pt, line width = 1pt, draw=black, fill=white]
\tikzstyle{vertexWj}=[circle,inner sep=0, minimum size =6 pt, line width = 1pt, draw=black!30, dashed, fill=white]

\tikzstyle{arete}=[draw,line width=1.5pt]
\tikzstyle{aretej}=[draw,dashed,color=black!30,line width=1.5pt]
\tikzstyle{aretegametree}=[draw,->, line width=1.2pt]

\tikzstyle{inv}=[circle,inner sep=0, minimum size =6 pt, line width = 1pt, draw=white, fill=white]

\tikzstyle{vertexB}=[circle,inner sep=0, minimum size =6 pt, line width = 1pt, draw=black, fill=black, text= white]
\tikzstyle{vertexBj}=[circle,inner sep=0, minimum size =6 pt, line width = 1pt, draw=black!30, dashed, fill=black!30, text= white]

\tikzstyle{redcircle}=[circle,inner sep=0, minimum size =15 pt, line width = 1pt, draw=red]

\tikzstyle{redcircled}=[circle,inner sep=0, minimum size =15 pt, line width = 1pt, draw=red, dashed]

\tikzstyle{redsquare}=[rectangle,inner sep=0, minimum size =15 pt, line width = 1pt, draw=red]

\tikzstyle{bluecircle}=[circle,inner sep=0, minimum size =15 pt, line width = 1pt, draw=blue]

\tikzstyle{bluesquare}=[rectangle,inner sep=0, minimum size =15 pt, line width = 1pt, draw=blue]

\tikzstyle{greencircle}=[circle,inner sep=0, minimum size =15 pt, line width = 1pt, draw=green]

\tikzstyle{greensquare}=[rectangle,inner sep=0, minimum size =15 pt, line width = 1pt, draw=green]

\tikzstyle{greentriangle}=[triangle,inner sep=0, minimum size =15 pt, line width = 1pt, draw=green]

\tikzstyle{redtriangle}=[triangle,inner sep=0, minimum size =15 pt, line width = 1pt, draw=blue]

\tikzstyle{bluetriangle}=[regular polygon,regular polygon sides=3,inner sep=0, minimum size =25 pt, line width = 1pt, draw=blue]

\begin{document}

\maketitle

\begin{abstract}
    The game {\sc influence} is a scoring combinatorial game that has been introduced in 2020 by Duchene et al \cite{duchene}. It is a good representative of Milnor's universe of scoring games, i.e. games where it is never interesting for a player to miss his turn. New general results are first given for this universe, by transposing the notions of mean and temperature derived from non-scoring combinatorial games. Such results are then applied to {\sc influence} to refine the case of unions of segments started in~\cite{duchene}. The computational complexity of the score of the game is also solved and proved to be PSPACE-complete. We finally focus on some specific cases of {\sc influence} when the graph is bipartite, by giving explicit strategies and bounds on the optimal score on structures like grids, hypercubes or torus.
\end{abstract}

%


\section{Introduction}


{\sc influence} is a two player scoring combinatorial game introduced by Duchene et al in \cite{duchene}. It is played on a directed graph where each vertex is colored either black or white. Alternately, both players (called Left and Right) take a vertex of their color: if it is Left's turn, she takes a black vertex, remove it and all its successors. If it is Right's turn, he takes a white vertex, remove it and all its predecessors. The game ends when the graph is empty. Each player scores the total number of vertices that he/she removed. Naturally, the objective of each player is to maximize his/her score.\\

The work presented in \cite{duchene} is the first study of a particular combinatorial game within the recent scoring framework of Larsson et al. \cite{Larsson2017}. Indeed, if scoring combinatorial game theory has been introduced by Milnor \cite{Milnor1953} and Hanner \cite{hanner1959} in the 50's, their work was not followed up, partly because of the hardness of solving such games. This is only since the last decade that new results appeared on the topic, in particular with the introduction of general frameworks of resolution, for particular families (also called {\em universes}) of scoring games.\\

In \cite{duchene}, the authors show that {\sc influence} belongs to the so-called {\em Milnor's universe} of scoring games. Roughly speaking, the major property that derived from this universe is that a player has never interest to miss his/her turn. They also partially solve the game on particular classes of digraphs, called {\em segments}, that correspond to alternated black and white oriented paths, with arcs from black to white vertices.\\

In this paper, we extend the case of segments to more general families of bipartite graphs having the same properties. For that purpose, we consider a particular version of {\sc influence} that is played on bipartite digraphs, where each set of the bipartition corresponds to vertices of the same color. In addition, all the arcs are oriented from black to white vertices. According to these constraints, this version amounts to playing the game on an undirected bipartite graph (where black and white vertices form the two sets of the bipartition), and where each player takes a vertex of his/her color, removes it and all its neighbours. To simplify the notation, this particular version of {\sc influence} will be called {\sc bipartite influence}. \\

\begin{figure}[ht]
    \centering
\scalebox{0.7}{\begin{tikzpicture}
\draw[arete] (0,0) grid (4,4);
\foreach \I in {0,...,4}
\foreach \J in {0,...,4}
\node[vertexW] at (\I,\J) {};

\foreach \I in {0,2,4}
\foreach \J in {0,2,4}
\node[vertexB] at (\I,\J) {};

\foreach \I in {1,3}
\foreach \J in {1,3}
\node[vertexB] at (\I,\J) {};

\node[redcircle] at (2,2) {};
\node[redcircled] at (1,2) {};
\node[redcircled] at (2,1) {};
\node[redcircled] at (3,2) {};
\node[redcircled] at (2,3) {};

\begin{scope}[shift={(7,0)}]
\foreach \I in {0,4}
{\node[vertexB](2\I) at (2,\I) {};
\foreach \J in {0,2,4}
\node[vertexB](\I\J) at (\I,\J) {};
}

\foreach \I in {1,3}
{
\node[vertexW](0\I) at (0,\I) {};
\node[vertexW](4\I) at (4,\I) {};
\node[vertexW](\I0) at (\I,0) {};
\node[vertexW](\I4) at (\I,4) {};
\foreach \J in {1,3}
\node[vertexB](\I\J) at (\I,\J) {};
}

\draw[arete] (00)--(01)--(02)--(03)--(04)--(14)--(24)--(34)--(44)--(43)--(42)--(41)--(40)--(30)--(20)--(10)--(00);
\draw[arete] (01)--(11)--(10) (03)--(13)--(14) (30)--(31)--(41) (34)--(33)--(43);

\node[redcircle] at (1,0) {};
\node[redcircled] at (0,0) {};
\node[redcircled] at (2,0) {};
\node[redcircled] at (1,1) {};

\end{scope}

\begin{scope}[shift={(14,0)}]

\node[vertexB](40) at (4,0) {};
\node[vertexB](24) at (2,4) {};
\node[vertexB](44) at (4,4) {};
\node[vertexB](42) at (4,2) {};
\node[vertexB](02) at (0,2) {};
\node[vertexB](04) at (0,4) {};

\foreach \I in {1,3}
{
\node[vertexW](0\I) at (0,\I) {};
\node[vertexW](4\I) at (4,\I) {};
\node[vertexW](30) at (3,0) {};
\node[vertexW](\I4) at (\I,4) {};}

\node[vertexB](13) at (1,3) {};
\node[vertexB](33) at (3,3) {};
\node[vertexB](31) at (3,1) {};

\draw[arete] (01)--(02)--(03)--(04)--(14)--(24)--(34)--(44)--(43)--(42)--(41)--(40)--(30);
\draw[arete] (03)--(13)--(14) (30)--(31)--(41) (34)--(33)--(43);
\end{scope}
\end{tikzpicture}}
\caption{First moves of \bip\ on a $5\times5$ grid \label{fig:exempleintro}}
\end{figure}
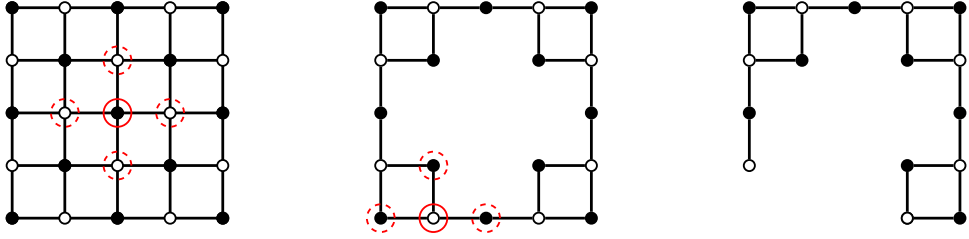

\bip\ is somehow very natural to play with, as bipartite bicolored graphs are a model for many 2-player abstract games. This is for example the case for checkerboard grids. Figure \ref{fig:exempleintro} illustrates the first two moves of a {\sc bipartite influence} game played on a $5\times 5$ checkerboard grid. After these first two moves, Left has scored 5 points (first move in the center), and Right 4 points (second move with the white bottom left vertex ). As discussed in \cite{duchene}, since isolated vertices can only be taken by the player who owns them, they are generally immediately removed from the graph and added to the score of the corresponding player.\\

In addition, this version of the game can be considered as a scoring variation of the well-known combinatorial game {\sc node-kayles} (in its partisan version). Indeed, in {\sc node-kayles}, each player removes a vertex and its neighborhood. The goal of each player is to be the last player to move. The correlation between {\sc node-kayles} and \bip\ is thus immediate.\\

For the rest of the paper, each instance of \bip\ will be denoted by a bipartite graph $G=(B\cup W,E)$, where $B$ and $W$ correspond respectively to the sets of black and white vertices, and every edge of $E$ is between a vertex of $B$ and $W$.\\
 
The main issue when studying a combinatorial game is the computation of the highest score of each player, by assuming both are playing optimally. This score is well defined as the game is finite and with no chance nor hidden information. The main results of this paper are described below. \\

In Section \ref{sec:material}, we provide the necessary material from scoring combinatorial theory that applies for the game. In particular and for the first time, we provide a formal framework inspired from Milnor and Hanner, and also from combinatorial (non-scoring) game theory \cite{aaron2013} for the concept of temperature and mean of a scoring game. 
In Section \ref{sec:influence}, we give general basic or known results on \bip. 
In Section \ref{sec:pspace}, we prove that {\sc influence} and even {\sc bipartite influence} are PSPACE-complete. The next section deals with {\sc bipartite influence} on graphs having some symmetry properties that induce a draw strategy. Section~\ref{sec:segments} is about a refinement of the results of \cite{duchene} on segments, in particular by using the concepts of mean and temperature. Finally, the game is considered on grids with two and three rows, with an exact value or a bound for the score in almost all cases. A conjecture is given for general grids.

\section{Background for scoring combinatorial game theory}\label{sec:material}

If scoring combinatorial game theory has been introduced in the 50's, the results of Milnor and Hanner have mainly been a bootstrap for combinatorial and economic game theory. A couple of recent studies (see~\cite{Larsson2017} for a survey) have settled general results for scoring game theory. For the reasons detailed in \cite{Larsson2017}, the one of Larsson et al. in \cite{Larsson2015} unifies the different previous frameworks. It has also more advantages, in particular because of its proximity with classical Combinatorial Game Theory. Our results will be presented according to this framework.

\subsection{Main definitions and notations}

The definitions below are derived from \cite{Larsson2015}. \\

A scoring game $G$ is recursively defined as a couple $\langle G^L|G^R\rangle$ where $G^L$ (resp. $G^R$) are two sets of scoring games. $G^L$ (resp. $G^R$) is called the set of Left (resp. Right) {\em options}. In terms of play, the Left (resp. Right) options of $G$ correspond to the games that can be reached by Left (resp. Right) in one move. 

If $G$ is a game with no Left option, then we write $G^L$ as $\emptyset^s$ (with $s$ being a real number) to indicate that if it is Left's turn, the game is over and the score of the game is $s$. The same notation applies for Right. Final positions correspond to games of the form  $\langle \emptyset^s|\emptyset^t \rangle$.\\

Note that in the rest of the paper, we will consider only {\em short} games, i.e. games $\langle G^L|G^R\rangle$ where the sets $G^L$ and $G^R$ are finite sets of short games.\\

In many games having symmetric rules, there exist no position of the form $\langle \emptyset^s|\emptyset^t \rangle$ with $s\neq t$. In particular, this is the case of {\sc influence}. As they correspond to final positions, games of the form $\langle \emptyset^s|\emptyset^s \rangle$ are called \em numbers \em and will be simply denoted by $s$. \\

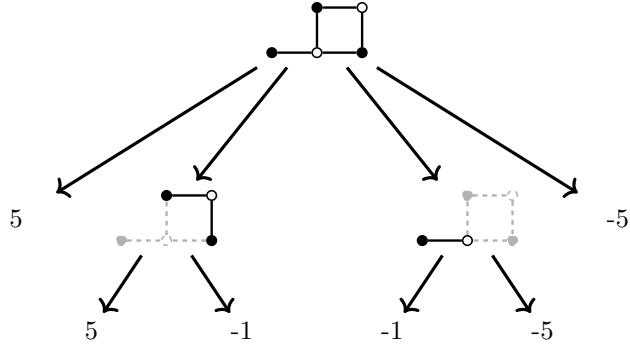
\begin{figure}
    \centering
    
\begin{tikzpicture}

\node(S) at (0,0) {
\scalebox{0.6}{\begin{tikzpicture}

\node[vertexB](0)  at (0,0) {};
\node[vertexW](1)  at (1,0) {};
\node[vertexB](2)  at (1,1) {};
\node[vertexW](3)  at (2,1) {};
\node[vertexB](4)  at (2,0) {};

\draw[arete] (0)--(1)--(2)--(3)--(4)--(1);
\end{tikzpicture}}
};

\node(1B) at (-2,-2.5) {
\scalebox{0.6}{\begin{tikzpicture}

\node[vertexBj](0)  at (0,0) {};
\node[vertexWj](1)  at (1,0) {};
\node[vertexB](2)  at (1,1) {};
\node[vertexW](3)  at (2,1) {};
\node[vertexB](4)  at (2,0) {};

\draw[aretej] (0)--(1)--(2) (4)--(1);
\draw[arete] (2)--(3)--(4);
\end{tikzpicture}}
};

\node(1W) at (2,-2.5) {
\scalebox{0.6}{\begin{tikzpicture}

\node[vertexB](0)  at (0,0) {};
\node[vertexW](1)  at (1,0) {};
\node[vertexBj](2)  at (1,1) {};
\node[vertexWj](3)  at (2,1) {};
\node[vertexBj](4)  at (2,0) {};

\draw[arete] (0)--(1);
\draw[aretej] (1)--(2)--(3)--(4)--(1);
\end{tikzpicture}}
};

\node[minimum size =30](1fB) at (-4,-2.5) {5};
\node[minimum size=30](1fW) at (4,-2.5) {-5};

\node(20) at (-3,-4) {5};
\node(21) at (-1,-4) {-1};
\node(22) at (1,-4) {-1};
\node(23) at (3,-4) {-5};

\path[aretegametree] (S) to (1B);
\path[aretegametree] (S) to (1W);
\path[aretegametree] (S) to (1fB);
\path[aretegametree] (S) to (1fW);
\path[aretegametree] (1B) to (20);
\path[aretegametree] (1B) to (21);
\path[aretegametree] (1W) to (22);
\path[aretegametree] (1W) to (23);

\end{tikzpicture}
    \caption{Example of a game tree from a position of \bip}
    \label{treeex}
\end{figure}

In scoring combinatorial game theory, the so-called {\em game tree} is a nice way to represent a game and all the possible moves. Given a game $G$, the set of nodes of its game tree are all the positions that can be reached when playing from $G$. The root of the tree is the starting position of the game. Given a node, its left (resp. right) sons are the Left (resp. Right) options. An example of game tree for \bip\ is given in Figure~\ref{treeex}. The scores on the leaves correspond to the difference between the number of vertices removed by Left and Right during the play.

If the notion of score is well defined when $G^L$ or $G^R$ is empty, Larsson et al.~\cite{Larsson2015} extended it recursively to any game $G$. Inspired by the notions of Left and Right stops in combinatorial game theory~\cite{aaron2013}, they defined in a similar way the notions of Left and Right scores corresponding to the score of the game depending on whether Left or Right starts the game.

\begin{definition}[Larsson et al \cite{Larsson2015}]
Let $G = \langle G^L | G^R \rangle$ be a scoring game. We call Left score and Right score of $G$, respectively denoted by $Ls(G)$ and $Rs(G)$, the following functions :

\begin{align*}
    &Ls(\langle \emptyset^s | G^R \rangle) = Rs(\langle G^L | \emptyset^s \rangle) = s \mbox{ if $s$ is a real number} \\
    &Ls(\langle G^L | G^R \rangle) = \underset{G^l \in G^L}{\max} Rs(G^l) \mbox{ if } G^L \neq \emptyset^s\\
    &Rs(\langle G^L | G^R \rangle) = \underset{G^r \in G^R}{\min} Ls(G^r) \mbox{ if } G^R \neq \emptyset^s\\
\end{align*}
\end{definition}

These definitions yield a partial order on the games. The convention is that positive scores are in favor of Left, whereas negative scores are for Right. Games with a score equal to zero are draws. On the example of Figure~\ref{treeex}, if $G$ is the root, we have $Ls(G)=5$ and $Rs(G)=-5$. This means that the player who starts wins the game.

The following definition concerns the length of a game. 

\begin{definition}
Let $G = \langle G^L | G^R \rangle$ be a scoring game. We define the length of $G$, denoted by $l(G)$ as the size of the longest sequence of moves in $G$.
\end{definition}

Recursively, $l(G)$ is thus defined as follows:
\begin{align*}
    l(\langle \emptyset^s|\emptyset^t\rangle) &= 0 \mbox{ for any real numbers $s$ and $t$} \\
    l(\langle G^L|G^R\rangle ) &= 1 + \underset{G' \in G^L\cup G^R}{ \max} l(G')
\end{align*}

As for standard combinatorial games, the disjunctive sum operator $+$ applied to scoring games is defined as follows:

\begin{definition}[Larsson et al. \cite{Larsson2015}]
Given two games $G_1$ and $G_2$, their \em disjunctive sum \em, written $G_1+G_2$ is defined as the following game:

\begin{align*}
 G_1+G_2 =   &\langle \emptyset^{\ell_1+\ell_2}|\emptyset^{r_1+r_2}\rangle \mbox{ if $G_1 = \langle \emptyset^{\ell_1}|\emptyset^{r_1}\rangle$ and $G_2 =  \langle \emptyset^{\ell_2}|\emptyset^{r_2}\rangle$} \\
   = &\langle \emptyset^{\ell_1+\ell_2}|G_1^R+G_2,G_1+G_2^R\rangle \mbox{ if $G_1 = \langle \emptyset^{\ell_1}|G_1^R\rangle$ and $G_2 = \langle \emptyset^{\ell_2}|G_2^R\rangle$ and at least $G_1^R$ or $G_2^R$ is not empty} \\
   = &\langle G_1^L+G_2,G_1+G_2^L|\emptyset^{r_1+r_2}\rangle \mbox{ if $G_1 = \langle G_1^L|\emptyset^{r_1}\rangle$ and $G_2 = \langle G_2^L|\emptyset^{r_2}\rangle$ and at least $G_1^L$ or $G_2^L$ is not empty} \\  
  =  &\langle G_1^L+G_2,G_1+G_2^L|G_1^R+G_2,G_1+G_2^R\rangle \mbox{ otherwise} \\
\end{align*}
\end{definition}

For the accuracy of this definition, note that $G_1^L+G_2$ does not exist if $G_1^L$ is empty. \\

In addition, the {\em opposite} (or {\em negative}) of a game $G$, denoted by $-G$, is defined as the game where the roles of Left and Right are exchanged. In \bip, it consists in exchanging the colors of the vertices. By definition of the scores, we have (\cite{duchene}, after Corollary 15) that for all $G$,
$$
Ls(-G)=-Rs(G).
$$


Finally, the notion of equivalence has been defined in order to replace big games by smaller ones in disjunctive sums, without changing the score:

\begin{definition}[Milnor \cite{Milnor1953}]
Two games $G_1$ and $G_2$ are \em equivalent \em (write $G_1 = G_2$) if for any game $G$, we have $Ls(G + G_1) = Ls(G + G_2)$ and $Rs(G+G_1) = Rs(G+G_2)$.
\label{def8}
\end{definition}

Note that this definition makes sense when $G$, $G_1$ and $G_2$ have similar properties (same universe - see next paragraph - or same ruleset).

\subsection{Milnor's Universe}

In scoring combinatorial game theory, a {\em universe} is generally defined as a set of games closed under disjunctive sum, by taking options and negative \cite{Larsson2015}. Several universes of scoring games have been considered in the literature (Stewart \cite{stewart2012}, Ettinger~\cite{ettinger}, Larsson \cite{Larsson2015}). In his paper, although it is not formulated in this way, Milnor introduced the universe of dicotic and nonzugzwang games \cite{Milnor1953}. Note that in \cite{Larsson2015}, this universe is also denoted $\mathbb{P}\mathbb{S}$. It is defined as follows:

\begin{definition}[Milnor \cite{Milnor1953}]
A game $G$ is \em dicotic \em if both players can move from every nonempty position in the game tree of $G$. 
\end{definition}

Formally, it is equivalent to say that for all $G = \langle G^L|G^R \rangle, G^L = \emptyset \Leftrightarrow G^R = \emptyset$.

\begin{definition}[Milnor \cite{Milnor1953}]
A game is \em nonzugzwang \em, or with \em no zugzwang \em, if for each position $G$ in its game tree, we have $Ls(G) \ge Rs(G)$.
\end{definition}

Roughly speaking, in a nonzugzwang game, a player never has interest in skipping his turn. \\

The disjunctive sum of games is natural for games that split into several small components during the play. Although the score of a sum of games is in general not equal to the sum of the scores of each game, it can be bounded according to the following theorem for the games that belong to Milnor's universe:

\begin{theorem}[Milnor \cite{Milnor1953}]
\label{thm:sum}
Let $G$ and $H$ be two dicotic nonzugzwang games,
we have 
$$
Rs(G) + Rs(H) \leq  Rs(G+H) \leq Ls(G)+Rs(H)\leq Ls(G+H)\leq Ls(G)+Ls(H).
$$
In particular, if $H$ is a number $s$ we have $Ls(G+s)=Ls(G)+s$ and $Rs(G+s)=Rs(G)+s$.
\end{theorem}

Being in Milnor's universe induces a couple of nice properties concerning the equivalence of games. The first result below means that every game $G$ of Milnor's universe has an inverse (i.e. $-G$), leading to the fact that this universe is an abelian group. This is not the case of the other universes studied in the literature. The second result yields a simple test to show that two games are equivalent.\\

\begin{lemma}\label{lemma equiv}
For any games $G$ and $H$ that are dicotic nonzugzwang, we have:
\begin{itemize}
\item $Ls(G-G)=Rs(G-G)=0$;
\item $Ls(G-H)=Rs(G-H)=0$ if and only if $G$ and $H$ are equivalent. 
\end{itemize}
\end{lemma}

As for standard combinatorial games, a relation order can be set on general scoring games (not necessary from Milnor's universe) \cite{ettinger}. In Milnor's universe, this relation can be defined as follows:

\begin{definition}
Let $G$ and $H$ be two games that are dicotic and nonzugzwang. We says that $G$ {\em dominates} $H$, and denote it by $G\geq H$, if $Rs(G-H)\geq 0$.
\end{definition}

If a game has dominated options, then it can be simplified by simply removing them:

\begin{theorem}\label{dominated moves}
Let $G=\langle G^L |G^R\rangle$ be a dicotic nonzugzwang game. Assume that there exist two options $G_1$ and $G_2$ in $G^L$ such that $G^l_2\geq G^l_1$. Then $G=\langle G^L\setminus G^l_1 | G^R\rangle$.

Similarly, if there exist $G^r_1\geq G^r_2$ in $G^R$, then $G=\langle G^L | G^R \setminus G^r_1 \rangle$.
\end{theorem}

\begin{proof}
By symmetry, we just prove the first part of the theorem. Let $G'=\langle G^L\setminus G^l_1 | G^R\rangle$.

We first prove that $G'$ belongs to Milnor's Universe. Since $G$ is dicotic and $G$ has at least two left options, $G'$ is still dicotic. All the positions of $G'$ are positions of $G$ except $G'$ itself. Thus to prove that $G'$ is nonzugzwang, we just need to prove that $Ls(G')\geq Rs(G')$.
Since $G^l_2$ dominates $G^l_1$, we have in particular $Rs(G^l_2)\geq Rs(G^l_1)$. Indeed, by Theorem \ref{thm:sum}, $Rs(G^l_2)\geq Rs(G^l_2-G^l_1)+Rs(G^l_1)$ and $Rs(G^l_2- G^l_1)\geq 0$.
Thus $\max_{G^l\in G^L}(Rs(G^l))=\max_{G^l\in G^L\setminus G^l_1} (Rs(G^l))$ and $Ls(G)=Ls(G')$.
On the other side, $Rs(G')=Rs(G)$ since the Right options of $G$ and $G'$ are the same. Thus $Ls(G')\geq Rs(G')$ and $G'$ is dicotic nonzugzwang.

Therefore, to prove that $G=G'$, one need to prove that $Ls(G-G')=Rs(G-G')=0$. Since $G-G'$ is dicotic nonzugzwang, we already have $Ls(G-G')\geq Rs(G-G')$. Thus one just need to prove that $Ls(G-G')\leq 0$ and that $Rs(G-G')\geq 0$.

Left can answer to any Right move in $G-G'$ by a move leading to a game $H-H=0$ since all the Right options in $G-G'$ have a symmetric Left option. Hence Left has a strategy to insure a score $0$ when playing second, which means that $Rs(G-G')\geq 0$.

Assume now that Left starts in the game $G-G'$. If Left does not play to $G_1^l-G'$, then Right can answer to a game $H-H=0$ and insure a score of $0$.
If Left play to $G_1^l-G'$, Right can answer to $G_1^l-G_2^l$ and insure a non positive score since $Ls(G_1^l-G_2^l)=-Rs(G_2^l-G_1^l)\leq 0$.
Thus Right has a strategy in second to insure non-positive score, which means that $Ls(G-G')\leq 0$.

\end{proof}

{\bf Application to \bip}\\

In \cite{duchene}, the authors proved that {\sc influence} belongs to Milnor's universe. It also naturally holds for \bip:

\begin{theorem}[\cite{duchene}]\label{thm:bipismilnor}
For any instance $G$ of \bip, we have that $G$ is dicotic and nonzugzwang.
\end{theorem}

Consequently, all the properties of Milnor's universe detailed above apply to \bip. In addition, Example \ref{exemple P5} below illustrates equivalent games in \bip, as well as dominated options. 

\begin{figure}[h]
    \centering
  \begin{tikzpicture}
\node[vertexB](1) at (0,0) {};
\draw (1) ++(0,0.3) node {$u$};
\node[vertexW](2) at (0.5,-0.5) {};
\node[vertexB](3) at (1,0) {};
\draw (3) ++(0,0.3) node {$v$};
\node[vertexW](4) at (1.5,-0.5) {};
\node[vertexB](5) at (2,0) {};

\draw[arete] (1)--(2)--(3)--(4)--(5);

\node at (1,-1) {$S_5$};

\node at (3,-0.25) {\huge $+$};

\begin{scope}[shift={(4,0)}]
\node[vertexB](1) at (0,0) {};
\node[vertexW](2) at (0.5,-0.5) {};
\node[vertexB](3) at (1,0) {};
\node[vertexW](4) at (1.5,-0.5) {};
\node[vertexB](5) at (2,0) {};

\node at (1,-1) {$S_5$};

\draw[arete] (1)--(2)--(3)--(4)--(5);
\end{scope}

\node at (7,-0.25) {\huge $=$};
\node at (8,-0.25) {\huge $2$};
\node at (9,-0.25) {\huge $+$};

\node[vertexB](6) at (10,0) {};
\node[vertexW](7) at (10.5,-0.5) {};
\draw[arete] (6)--(7);

\node at (10.25,-1) {$S_2$};
 
    \end{tikzpicture}
    \caption{In \bip, the game $S_5+S_5$ is equivalent to the game $2+S_2$.}
    \label{fig:exP5}
\end{figure}
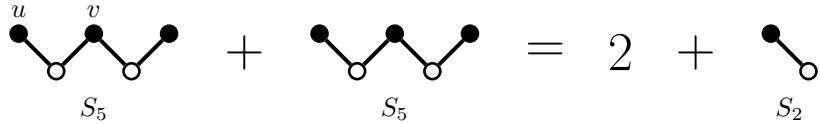

\begin{exemple} 
\label{exemple P5}
Denote by $S_5$ the path of length 5, with the two extremities and the middle vertex colored black, and the two other vertices colored white. Denote also by $S_2$ the path of length two with one black and one white vertex. See Figure \ref{fig:exP5}.

We will prove that $2S_5 = 2 + S_2$. According to Lemma~\ref{lemma equiv}, it is sufficient to prove $Ls(2S_5 - (2 + S_2)) = Rs(2S_5 -(2+S_2)) = 0$. As $S_2 = - S_2$ and as numbers can be removed from a sum according to Theorem~\ref{thm:sum}, we will prove $Ls(2S_5 + S_2) = Rs(2S_5 + S_2)  = 2$.
As a first step, when considering the Left options from $2S_5+S_2$, playing an extremity of a $S_5$ (i.e. the vertex $u$ in Figure~\ref{fig:exP5}) always leads to a position that is dominated by playing the middle of $S_5$ (i.e. the vertex $v$).  Therefore, when computing $Ls(2S_5+S_2)$, one only take into account the option $v$ (that yields $5$ points) and the move on $S_2$ (that yields $2$ points):

\begin{align*}
    Ls(2S_5 + S_2) &= \max \left ( 5+ Rs(S_5 + S_2), 2 + Rs(2S_5) \right ) \\
    &= \max \left ( 5+ \min \big( -2 + Ls(S_5), -3 + Ls(S_2+S_2) \big), 2 - 3 + Ls(S_5 + S_2) \right ) \\
    &= \max \left ( 5+ \min \big( 3, -3 \big), -1 + \max \big( 5 + Rs(S_2), 2 + Rs(S_5) \big)  \right ) \\
   &= \max \left ( 2, -1 + \max \big( 3, 2 - 3 + Ls(S_2) \big)  \right ) \\
    &= \max \left ( 2, -1 + \max \big( 3, - 1 + 2 \big)  \right ) \\
    &= 2. \\
  & \\
    Rs(2S_5 + S_2) &= \min \left ( -3+ Ls(S_5 + S_2 + S_2), -2 + Ls(2S_5) \right ) \\
     &= \min \left ( -3+ Ls(S_5), -2 + 5 + Rs(S_5) \right ) \\
     &= \min \left ( 2, 3 - 3 + Ls(S_2) \right ) \\
    &= \min \left (2, 2\right ) \\
     &= 2.
\end{align*}

This proves that $2S_5 = 2 + S_2$. Moreover, as $S_2 = -S_2$, we also have that $4S_5 = 2 + S_2 + 2 + S_2 = 4$.
\end{exemple}

\subsection{Mean and Temperature}

In this section, we develop the so-called concepts of {\em mean} and {\em temperature} of a game. These two notions have been introduced by Milnor and Hanner~\cite{Milnor1953, hanner1959} for scoring games, and then have been formalized in the context of combinatorial game theory by Conway~\cite{conway}. In his book, Siegel~\cite{aaron2013} gives formal definitions and proofs about them. In this section, we propose to do a similar work in the context of scoring combinatorial game theory. Our results cover the results of Hanner and extend them by following the CGT formalism. In addition, some of them will be used in Section~\ref{sec:segments} to yield results about \bip\ on segments.


While playing on a (disjunctive) sum of games, general information can be deduced from the information on each term of the sum. As said previously, the Left score of a sum of games is not generally the sum of the Left scores. Nevertheless, Milnor and Hanner have introduced several tools to get some knowledge on sums of games. The most intuitive one is the \em mean \em of a game, which is very useful if the same game is repeated many times in a sum.

\begin{definition}[Hanner \cite{hanner1959}]
The \em left (resp. right) mean \em of a game $G$, denoted $m_L$ (resp. $m_R$), is the score obtained in average by playing $G$ a huge number of times if Left (resp. Right) starts, if it exists. More formally:

\begin{align*}
    m_L(G) &= \underset{n \to \infty}{\lim} \frac{Ls(nG)}{n} \\
    m_R(G) &= \underset{n \to \infty}{\lim} \frac{Rs(nG)}{n}
\end{align*}
where $nG$ corresponds to the disjunctive sum $G+\cdots G$ with $n$ terms.
\label{def:mean}
\end{definition}

In Milnor's universe, both left and right means always exist:

\begin{theorem}[Milnor \cite{Milnor1953}]
Let $G$ be a dicotic nonzugzwang game. The values $m_L(G)$ and  $m_R(G)$ exist and are equal.
\label{thmmoy}
\end{theorem}

The proof derives from the fact that asymptotically, the first move is made on a game that cannot significantly modify the rest of the sum. According to this result, in Milnor's universe, the mean value is thus unique and will be denoted $m(G)$.\\

The mean of a game can be interpreted as what a player can expect from a game in a large sum. However, it gives no clue about which game of the sum it is strategic to play in. In order to answer to this question, Hanner~\cite{hanner1959} introduced the so-called parameter $\si(G)$. The more $\si(G)$ is high, the  more the game is said \em hot \em, meaning a player will want to play in it in priority. This parameter is similar to the \em temperature \em notion known in standard combinatorial games \cite{aaron2013}. Nevertheless, Hanner's definition of $\si$ requires the use of an algorithm and is not very tractable. In what follows, we propose a formalization of this notion of temperature for scoring games, by adapting the formalism used in combinatorial games. Many results from Hanner \cite{hanner1959} about $\sigma$ remain true within this formalism.

\begin{definition}\label{def:temp}
Let $G = \langle G^L|G^R\rangle$ be a dicotic nonzugzwang game. Given a real number $t\geq 0$, we define $G$ {\em cooled} by $t$, denoted by $G_t$, as follows. If $G$ is a real number $k$, put $G_t = k$ and $t_0=0$. Otherwise, let
$\widetilde{G}_t=\langle G_t^L - t|G^R_t + t \rangle$,
$t_0=\mbox{min } \{ t \ge 0 | Ls(\widetilde{G}_t)= Rs(\widetilde{G}_t)\}$
and define $$G_t = \left\{
    \begin{array}{ll}
         &  \widetilde{G}_t \mbox{ if } t\le t_0\\
         & Ls(G_{t_0}) \mbox{ otherwise} 
    \end{array} 
    \right.$$
    
The value $t_0$ is called the {\em temperature} of $G$ and is denoted by $\si(G)$.
\end{definition}

If $t$ is higher or equal to the temperature of $G$, the game $G_t$ is said to be \textit{frozen} as it behaves as a number. If $t$ is lower than the temperature of $G$, playing in $G_t$ can be understood as "playing in $G$ but each move costs $t$ to the player". \\

This definition of temperature matches with the usual definition of temperature in combinatorial games without scores. In both cases, the temperature is defined as a minimum value, which may not exist.

\begin{figure}[h]
    \centering
    
\begin{tikzpicture}

\begin{scope}[shift={(-3.5,0)}]
\node[vertexB](1) at (0,0) {};
\node[minimum size=20](2) at (-1,-1.5) {5};
\node[vertexB](3) at (1,-1.3) {};
\node[minimum size=20](4) at (0.5,-2.75) {-1};
\node[minimum size=20](5) at (1.5,-2.75) {-5};

\path[aretegametree] (1) to (2);
\path[aretegametree] (1) to (3);
\path[aretegametree] (3) to (4);
\path[aretegametree] (3) to (5);

\node at (0,-4) {Tree of $S_5$};
\end{scope}

\begin{scope}[shift={(2.5,0)}]
\node[vertexB](1) at (0,0) {};
\node[minimum size=20](2) at (-1,-1.5) {$5-t$};
\node[vertexB](3) at (1,-1.3) {};
\node[minimum size=20](4) at (0.3,-2.75) {-1};
\node[minimum size=20](5) at (1.75,-2.75) {$-5+2t$};

\path[aretegametree] (1) to (2);
\path[aretegametree] (1) to (3);
\path[aretegametree] (3) to (4);
\path[aretegametree] (3) to (5);

\node[text centered, text width=3cm] at (0.3,-4) {Tree of $S_5$ cooled by $t$ for $t\in [0,2]$};
\end{scope}

\begin{scope}[shift={(8.5,0)}]
\node[vertexB](1) at (0,0) {};
\node[minimum size=10](2) at (-1,-1.5) {$5-t$};
\node[minimum size=10](6) at (1,-1.5) {$-3+t$};

\path[aretegametree] (1) to (2);
\path[aretegametree] (1) to (6);

\node[text centered, text width=3cm] at (0.3,-4) {Tree of $S_5$ cooled by $t$ for $t\in [2,4]$};
\end{scope}
\end{tikzpicture}
\caption{Tree of the game $S_5$ and the same game cooled by $t$ for different values of $t$.}
\label{fig:treecooled}
\end{figure}
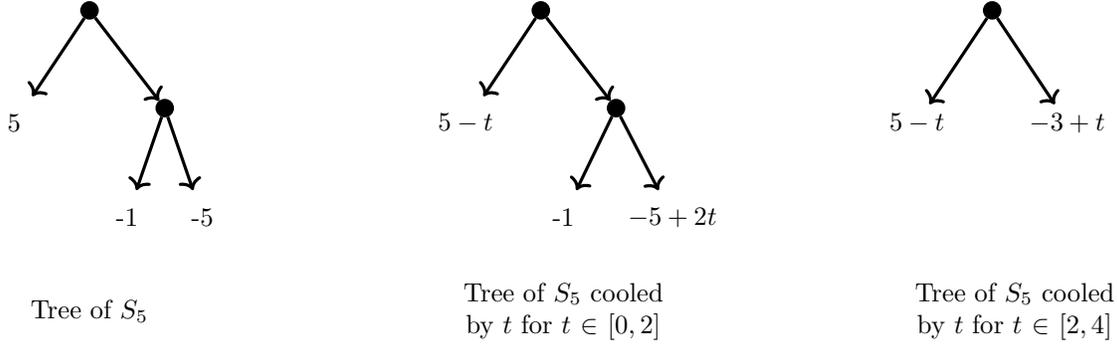

Before proving that the temperature is well-defined, we give an example of the computation of the temperature.
\begin{exemple}\label{ex:temperatureS5}
We compute in this example the temperature of the instance $S_5$ of \bip\ defined in Example \ref{exemple P5}. The tree of $S_5$ is given in Figure~\ref{fig:treecooled} where the dominated moves corresponding to the extremities have been removed. We have $S_5=\langle 5|\langle -1|5\rangle \rangle$.

Let $t\geq 0$. By definition of a cooled game, we have
$$
\widetilde{(S_5)}_t=\langle 5_t - t |\langle -1, -5\rangle_t+t\rangle.$$

Since numbers are already cooled, we have
$$
\widetilde{(S_5)}_t=\langle 5 - t | \langle -1, -5\rangle_t+t\rangle.
$$
The temperature of $\langle -1|-5\rangle$ is 2. Hence, for $t\leq 2$, we have that $\langle -1|-5\rangle_t + t=\langle -1|-5+2t\rangle$, and $\widetilde{(S_5)}_t=\langle 5 - t | \langle -1|-5+2t\rangle \rangle$ which corresponds to the second case of Figure~\ref{fig:treecooled}. For $t\geq 2$, the game $\langle -1|-5\rangle_t$ is frozen and equals to $-3$. Thus $\widetilde{(S_5)}_t=\langle 5 - t | -3+t\rangle$ which corresponds to the rightmost case of Figure~\ref{fig:treecooled}. When $t=4$, we have for the first time $Ls(\widetilde{(S_5)}_t)=Rs(\widetilde{(S_5)}_t)$ and thus the temperature of $(S_5)$ equals $4$. After $t=4$, $(S_5)_t$ is frozen and equals to $1$. 

\end{exemple}

\begin{proof}[Proof of validity of Definition~\ref{def:temp}]
We prove by induction on the length of $G$ that for any $t \ge 0$, $G_t$ exists and the two maps $t \mapsto Rs(G_t)$ and $t \mapsto Ls(G_t)$ are continuous.  

If $G$ is a number, its temperature is 0 and is well defined. Hence $G_t$ exists for all $t \ge 0$ and the two score functions are constant.

Otherwise, $G = \langle G^L|G^R\rangle $ with $G^L$ and $G^R$ non empty since the game is dicotic. By induction, the temperature of all the games in $G^L$ and $G^R$ are well defined. Thus the games in $G^L_t$ and $G^R_t$ exist for all $t \ge 0$ and the corresponding score functions are continuous.

Therefore, the functions $t \mapsto Ls(\widetilde{G_t})$ and $t \mapsto Rs(\widetilde{G_t})$ exist and are continuous as maximum and minimum of (a finite number of) continuous functions. 

Consider now the function $f:t \mapsto Ls(\widetilde{G_t})-Rs(\widetilde{G_t})$ defined on $\R^+$. The function $f$ is well defined, continuous and 
$$f(t)=\underset{G^l \in G^L}{ \max } \big( Rs(G^l_t) \big) - \underset{G^r \in G^R}{ \min } \big( Ls(G^r_t) \big) - 2t.$$

We have $f(0) = Ls(G) - Rs(G) \ge 0$ as $G$ is nonzugzwang. Furthermore, $\underset{t \to \infty}{\mbox{lim }} f(t) = -\infty$ since the functions $t\mapsto Rs(G^l_t)$ and $t\mapsto Ls(G^r_t)$ are in a finite number and, by definition, constant for $t$ large enough.
Since $f$ is continuous, there exists a minimum $t_0$ such that $f(t_0)=0$  which guarantees the existence of the temperature and the game $G_t$.

Moreover, for any $0 \le t \le t_0$, $Ls(G_t)$ and $Rs(G_t)$ are continuous as minimum and maximum of continuous functions; for $t \geq t_0$, they are constant so continuous. Finally, $Ls(G_t)$ and $Rs(G_t)$ are continuous in $t_0$ since $Ls(G_{t_0}) = Rs(G_{t_0})$.
\end{proof}

With a similar induction, one can prove the following corollary.
\begin{corollaire}\label{cor:decreasing}
Let $G$ be dicotic nonzugzwang game. Then for any $t\geq 0$, $G_t$ is a dicotic nonzugzwang game.
The function $t\mapsto Ls(G_t)$ and $t\mapsto Ls(\widetilde{G_t})$ are decreasing functions.
Symmetrically, $t\mapsto Rs(G_t)$ and $t\mapsto Rs(\widetilde{G_t})$ are increasing functions.

In particular, if $t\geq \si(G)$, we have $Ls(\widetilde{G_t})\leq G_{\si(G)} \leq Rs(\widetilde{G_t})$.
\end{corollaire}

A natural question is about how a cooled game is far from the original game, in terms of Left and Right scores. The following property yields an answer.

\begin{proposition}
Let $G$ be a dicotic nonzugzwang game and let $t \ge 0$. We have $ 0 \le Ls(G) - Ls(G_t) \le t$ and $0 \ge Rs(G) - Rs(G_t) \ge -t$

\label{propcool}
\end{proposition}

\begin{proof}
By considering $-G$ instead of $G$, it suffices to prove the result for $Ls$ since $Rs(-G)=-Ls(G)$ and $Rs(-G_t)=-Ls(G_t)$.
By Corollary \ref{cor:decreasing}, $t\mapsto Ls(G_t)$ is decreasing, which implies that $Ls(G)-Ls(G_t)\geq 0$.
It remains to prove that $Ls(G)-Ls(G_t)\leq t$

If $G$ is a number, we have for any $t \ge 0$, $G_t = G$, which implies the result.

Otherwise, consider the case $0\le  t\le \si(G)$. 
Let $G^{l}$ be an optimal move for Left in $G$, we have $Ls(G)=Rs(G^l)$. Since $Ls(G_t)=\max_{l'\in L}(Rs(G^{l'}_t)-t)$, we have in particular $Ls(G_t)\geq Rs(G^l_t)-t$.
Thus, $Ls(G)-Ls(G_t)\leq Rs(G^l)-Rs(G^l_t)+t$. Since $t\mapsto Rs(G^l_t)$ is increasing, $Rs(G^l)-Rs(G^l_t)\leq 0$ leading to $Ls(G)-Ls(G_t)\leq t$.

Consider now the case $t > \si(G)$, we have $Ls(G) - Ls(G_t) = Ls(G) - Ls(G_{\si(G)})$. By the last inequality above, we have $ t > \si(G) \ge Ls(G) - Ls(G_t)$.



\end{proof}



As said before, this definition of a cooled game can be understood as a game where moving has a cost. Hence, the temperature is defined as the highest cost that can be set such that the game remains nonzugzwang. Since the temperature has been defined to handle sums of games, we state preliminary results  about the sum of two cooled games, starting by the case where one the games is a number.


\begin{lemma}\label{lem:linearnumber}
Let $G$ be a dicotic nonzugzwang game. Let $s$ be a number. Then $(G+s)_t=G_t+s$. In particular, $\si(G+s)=\si(G)$.
\end{lemma}

\begin{proof}
The result is prove by induction on the length of $G$. If $G$ is a number, so is $G+s$ and we have $(G+s)_t=G+s=G_t+s$.

Otherwise:
\begin{align*}
   \widetilde{(G+s)}_t &= \langle (G^L+s)_t-t | (G^R+s)_t+t\rangle\\
    &=\langle G^L_t+s-t | G^R_t+s+t\rangle \mbox{ by induction }  \\
    &=\langle G^L-t| G^R+t\rangle +s \\
    &=\widetilde{G}_t +s.\\
\end{align*} 

In particular, $Ls(\widetilde{(G+s)_t})=Ls(\widetilde{G}_t)+s$ and 
$Rs(\widetilde{(G+s)_t})=Rs(\widetilde{G}_t)+s$. This implies that $\si(G+s)=\si(G)$ and that  $(G+s)_t=G_t+s$. 
\end{proof}

Despite the above result, it is not clear in our opinion that the cooling function is a homomorphism for any two scoring games. If this result is true in the CGT context, a direct adaptation in Milnor's universe is not guaranteed, as $\widetilde{G_t}$ can be outside this universe. However, if the cooled games are not equivalent, the result below show that it is true for the scores. In addition, as for CGT, we also show that the temperature is a submaximal function.

\begin{theorem}\label{thm:sumtemp}
Let $G$ and $H$ be two dicotic nonzugzwang games. We have $Ls((G + H)_t) =Ls(G_t + H_t)$ and $Rs((G+H)_t)=Rs(G_t+H_t)$. Moreover, $\si(G+H)\leq \max(\si(G),\si(H))$ and if $\si(H)<\si(G)$, then $\si(G+H)=\si(G)$. 
In particular, we have: $(G+H)_{\si(G+H)}=G_{\si(G)}+H_{\si(H)}$
\end{theorem}

\begin{proof}
We prove the proposition by induction on the length of $G+H$.

If $G$ or $H$ is a number, the proposition is a direct consequence of Lemma~\ref{lem:linearnumber}. 



From now on, we assume that both $G$ and $H$ are not numbers, and, without loss of generality, we assume that $\si(G) \ge \si(H)$.
We first prove that $Ls(\widetilde{(G+H)_t})=Ls(G_t+H_t)$ for any $t\leq \si(G)$.

Let $t \ge 0$.
\begin{align*}
    Ls(\widetilde{(G+H)}_t) &= \max_{G^l\in G^L,H^l\in H^L}\{Rs((G^l+H)_t-t), Rs((G+H^l)_t-t)\}\\
    &= \max_{G^l\in G^L,H^l\in H^L}\{Rs((G^l+H)_t), Rs((G+H^l)_t)\}-t\\
    &= \max_{G^l\in G^L,H^l\in H^L}\{Rs(G^l_t+H_t), Rs(G_t+H^l_t)\}-t \mbox{ by induction }  \\
    &= \max_{G^l\in G^L,H^l\in H^L}\{Rs(G^l_t+H_t-t), Rs(G_t+H^l_t-t)\}.\\
\end{align*} 

Assume first that $t\le \si(H)$, then
\begin{align*}
G_t + H_t &=\langle G^L_t-t|G^R_t+t\rangle + \langle H^L_t-t|H^R_t+t\rangle\\
Ls(G_t+H_t)&= \max_{G^l\in G^L,H^l\in H^L}\{Rs(G^l_t+H_t-t), Rs(G_t+H^l_t-t)\}
\end{align*} and  $Ls(\widetilde{(G+H)}_t)=Ls(G_t + H_t)$.

Assume now that $\si(H)\le t \leq \si(G)$. In particular, $H_t$ is a number equal to $H_{\si(H)}$. We prove that there is an option in $G^L_t+H_t-t$ that has a Right Score greater all the options in $G_t + H^L_t-t$. 
Let $G^l\in G^L$ such that $Ls(G_t)=Rs(G^l_t-t)$ ($G^l$ exists since $G_t$ is not a number). Let any $H^l\in H^L$.

Since $t\geq \si(H)$ and since $t\mapsto Ls(\widetilde{H_t})$ is, by Corollary \ref{cor:decreasing}, a decreasing function, we have $Rs(H^l_t-t)\leq Ls(\widetilde{H_t})\leq Ls(\widetilde{H_{\si(H)}})=H_{\si(H)}$.

Therefore, we have the following inequalities:
\begin{align*}
    Rs(G_t+H^l_t-t) &\leq Ls(G_t)+Rs(H^l_t-t) \mbox{ since all the components are in Milnor's universe {(see \ref{thm:sum} }}\\
    &\leq Rs(G^l_t-t) + H_{\si(H)}\\
    &\leq Rs(G^l_t-t+H_t) \mbox{ since } H_t=H_{\si(H)}\\
\end{align*} 

Thus, the Left score of $\widetilde{(G+H)}_t$ is computed only on the options $G^L+H_t-t$ and thus:
\begin{align*}
    Ls(\widetilde{(G+H)}_t)&=\max_{G^l\in G^L}\{Rs(G^l_t+H_t-t\}\\
    &= \max_{G^l\in G^L}\{Rs(G^l_t-t)\} +H_{\si(H)}\\
    &=Ls(G_t)+H_{\si(H)} \mbox{ since } t\leq\si(G)\\
    &=Ls(G_t+H_t).\\
\end{align*}

We have proved that if $t\leq \si(G)$, then $Ls(\widetilde{(G+H)_t})=Ls(G_t+H_t)$.
In a similar way, one can prove that $Rs(\widetilde{(G+H)_t})=Rs(G_t+H_t)$.

The two functions $t\mapsto Ls(\widetilde{(G+H)_t})$ and $t\mapsto Rs(\widetilde{(G+H)_t})$ are respectively decreasing and increasing functions. They are both equal to $G_{\si(G)}+H_{\si(H)}$ in $t=\si(G)$. Thus the first time they are equal is before (or in) $\si(G)$ which means precisely that $\si(G+H)\leq \si(G)$. Furthermore, we have:
\begin{itemize}
    \item for $0\leq t <\si(G+H)$, $Ls((G+H)_t)=Ls(\widetilde{(G+H)_t})=Ls(G_t+H_t)$ and
$Rs((G+H)_t)=Rs(\widetilde{(G+H)_t})=Rs(G_t+H_t)$;
\item for $\si(G+H)\leq t\leq \si(G)$, $t\mapsto Ls(\widetilde{(G+H)_t})$ and $t\mapsto Rs(\widetilde{(G+H)_t})$ are constant and both equal to the value in $\si(G+H)$ and $\si(G)$. In particular, $$Ls((G+H)_t)=Ls(\widetilde{(G+H)_t})=Ls(G_t+H_t)=(G+H)_{\si(G)+\si(H)}=G_{\si(G)}+H_{\si(H)},$$
$$Rs((G+H)_t)=Rs(\widetilde{(G+H)_t})=Rs(G_t+H_t)=(G+H)_{\si(G)+\si(H)}=G_{\si(G)}+H_{\si(H)}.$$
\end{itemize}

If $\si(H)<\si(G)$, then for any $\si(H)<t<\si(G)$, $$Ls(\widetilde{(G+H)_t})=Ls(G_t)+H_{\si(H)}>Rs(G_t)+H_{\si(H)}=Rs(\widetilde{(G+H)_t}),$$and therefore $\si(G+H)=\si(G)$.

Finally, for $t\leq \si(G+H)$,  $Ls((G+H)_t)=(G+H)_{\si(G)+\si(H)}$ and $Ls(G_t+H_t)=G_{\si(G)}+H_{\si(H)}$. From what precedes, the two values are equal, which concludes the proof.
\end{proof}







The score after having cooled a game is defined as the score obtained after having increased the cost of the plays. At the temperature, one can also consider it as the score on which Left and Right agree to play the game regardless who starts. In this sense, this definition is close to the mean of the game. As first announced in \cite{hanner1959} in another format, the following result formalizes the connection between the two notions.

\begin{proposition}
Let $G$ be a nonzugzwang dicotic game. We have $G_{\si(G)} = m(G)$. 
\label{prop:sim}
\end{proposition}

\begin{proof}
We denote by $x$ the number $G_{\si(G)}$. We prove that $m(G)=x$
By Theorem \ref{thm:sumtemp}, $\si(nG)\leq \si(G)$ and $(nG)_{\si(nG)}=nG_{\si(G)}=nx$.
In particular, for $t$ large enough, $(nG)_t=nx$.

Therefore, from Proposition~\ref{propcool}, we have $Ls(nG) - t \le Ls((nG)_t) = nx = Rs((nG)_t) \le Rs(nG) + t$

By dividing by $n$, we have $$\frac{Ls(nG)}{n} - \frac{t}{n} \le x \le \frac{Rs(nG)}{n} + \frac{t}{n}$$.

Finally, as $\underset{n \to \infty}{lim } \frac{Ls(nG)}{n} = \underset{n \to \infty}{lim } \frac{Rs(nG)}{n} = m(G)$, and $\frac{t}{n} \underset{n \to \infty}{\to} 0$, we have $x = m(G)$.
\end{proof}

Next theorem expresses in which sense a small temperature means that playing first is not very important.

\begin{theorem}\label{thm:tempmean}
Let $G$ be a nonzugzwang dicotic game, we have $m(G) - \si(G) \le Rs(G) \le m(G) \le Ls(G) \le m(G) + \si(G)$.
\end{theorem}

\begin{proof}
As $G$ is supposed to be dicotic and nonzugzwang, by symmetry, it suffices to prove $m(G) \le Ls(G) \le m(G) + \si(G)$. 
Indeed, by considering $-G$ instead of $G$, we have: 
\begin{align*}
&m(-G) \le Ls(-G) \le m(-G) + \si(-G) \\
\Leftrightarrow \hspace{2cm}& -m(G) \le -Rs(G) \le -m(G) + \si(G) \mbox{ as } \si(G) = \si(-G) \mbox{ and } m(G)=-m(G)\\
\Leftrightarrow \hspace{2cm}& m(G) \ge Rs(G) \ge m(G) - \si(G) 
\end{align*}

By Corollary \ref{cor:decreasing}, the function $t\mapsto Ls(G_t)$ is decreasing. Since $G_0=G$ and $Ls(G_t)=m(G)$ for $t$ large enough, we have $Ls(G)\geq m(G)$.

By Proposition \ref{propcool}, $Ls(G)-Ls(G_{\si(G)})\leq \si(G)$ which is equivalent to $Ls(G)\leq m(G)+\si(G)$.



\end{proof}

Theorem \ref{thm:tempmean} can be applied to sum of games to bound the overall score according to the mean and the temperature of each term of the sum. The following corollary confirms Theorem 1 of \cite{hanner1959} with our definition of the temperature. We will use this result for \bip\ on sum of segments in Section \ref{sec:segments}.

\begin{corollaire}
Let $G = G_1 + \dots + G_n$ be a sum of dicotic nonzugzwang games. Let $m_i$ be the mean of $G_i$ and $\si = \underset{1\le i \le n}{\max} \si(G_i)$. We have $$\underset{i=0}{\overset{n}{\sum}} m_i - \si \le Rs(G) \le \underset{i=0}{\overset{n}{\sum}} m_i \le Ls(G) \le \underset{i=0}{\overset{n}{\sum}} m_i + \si .$$
\label{coro:sumtemp}
\end{corollaire}

\begin{proof}

This result is straightforward as we know that $m(G) = \underset{i=1}{\overset{n}{\sum}} m(G_i)$ (as a sum of limits) and $\si(G) \le \underset{1\le i \le n}{\max} \si(G_i)$ (from Theorem~\ref{thm:sumtemp}).

\end{proof}

\section{General results on \bip}\label{sec:influence}

The current section gives general results about \bip\ that will be useful in the rest of the paper. We first recall a property introduced in \cite{duchene} about the so-called relevant graphs. In the context of \bip, the vertices from a graph $G$ that can be directly removed by playing a vertex $v$ are $v$ itself, its neighbors and the isolated vertices that possibly appeared (with the same color as $v$). We denote by $\rmv{G}{v}$ this set of removed vertices. Consequently, the score functions can also be expressed as follows:
$$
Ls(G)=  \max_{x \in B}\bigl\{ \vert \rmv{G}{x} \vert  +
Rs(G\setminus \rmv{G}{x})\bigl\} \quad \mbox{and} \quad Rs(G)=  \min_{y \in W} \bigl\{- \vert \rmv{G}{y} \vert   +
Ls(G\setminus \rmv{G}{y})\bigl\}.
$$

\subsection{Included moves}

 In {\sc Influence}, it has been proved that a move included in another is dominated by it \cite{duchene}. More formally, we have the following theorem that is a restatement of Corollary~13 in \cite{duchene}.

\begin{theorem}[\cite{duchene}]
\label{thm:includedmoves}
Let $G = (B \cup W,E)$ be a bipartite graph and $u,v$ two vertices of the same colour such that $Rmv(G,u) \subseteq Rmv(G,v)$. Then, the option of playing $v$ dominates the option of playins $u$. Otherwise said, playing $v$ is always at least as good as playing $u$. 


\end{theorem}

\subsection{Removing and adding vertices}

A way to obtain bounds on the score is to consider that some vertices are automatically given to the opponent, in order to simplify the structure of the graph without those vertices. The following result explains how the score may vary when sets of vertices of the same color are removed.

\begin{theorem}\label{thm:givingvertices}
Let $G = (B\cup W,E)$ be a bipartite graph. Let $B_0 \subseteq B$ and $W_0 \subseteq W$.

We have: \begin{align}
    Ls(G) &\le Ls(G\bs W_0) + |W_0| \\
    Rs(G) &\ge Rs(G\bs B_0) - |B_0| \\
    Ls(G) &\ge Ls(G\bs B_0) - |B_0| \\
    Rs(G) &\le Rs(G\bs W_0) + |W_0|.
\end{align}

\label{thmineq}
\end{theorem}

Roughly speaking, it can be explained by the fact that removing a set of white vertices and giving them to Left may prevent Right from doing some good moves but does not prevent Left to apply the same strategy. 


\begin{proof}
By symmetry, it suffices to prove $(1)$ and $(4)$. We do this by induction on $|V|$. The result is trivially true if $W_0$ is empty, if $|V|=1$, or if the first player has no move. Assume it is true for any graph with at most $n \ge 1$ vertices and where the first player has a move. 

Let $G = (B\cup W,E)$ be a graph with $n+1$ vertices and let $W_0 \subset W$ be a set of white vertices. The main idea of the proof is that by applying the optimal strategy existing in $G$ on the graph with the removed vertices, the player scores at least the same number of points. \\

We first prove $(1) :$ let $v \in B$ be a vertex played by Left in an optimal move in $G$. Denote by $W_1$ the white vertices of $W_0$ removed by playing $v$. Namely, we have $W_1 = W_0 \cap \rmv{G}{v}$, and define $W_2$ as the other white vertices of $W_0$, i.e. $W_2 = W_0 \bs W_1$. 

\noindent Consider now a strategy for Left that starts by playing $v$ in $G\bs W_0$. By definition of $Ls$, we have :

$$Ls(G\bs W_0) \ge Rs \big ( (G\bs W_0)\bs \rmv{G\bs W_0}{v}\big ) + \big |\rmv{G\bs W_0}{v} \big |$$

\noindent Since $W_1 = W_0 \cap \rmv{G}{v} \subset \rmv{G}{v}$, we have $W_0 \cup \rmv{G}{v} = W_1 \cup W_2 \cup \rmv{G}{v} = W_2 \cup \rmv{G}{v}$. Therefore, we have:

$$Ls(G\bs W_0) + |W_0| \ge |\rmv{G\bs W_0}{v}| + Rs\big ( (G\bs W_2)\bs \rmv{G}{v} \big ) + |W_1| + |W_2|$$

\noindent Now, by induction we know that $Rs(G\bs \rmv{G}{v}) \le Rs((G \bs \rmv{G}{v})\bs W_2) + |W_2|$. Thus, by applying this equality to the previous one, we get $  Ls(G\bs W_0) + |W_0| \ge |\rmv{G\bs W_0}{v}| + |W_1| + Rs(G\bs \rmv{G}{v} $.

\noindent Finally, by the triangular inequality $|\rmv{G\bs W_0}{v}| + |W_1| \le |\rmv{G\bs W_0}{v} \cup W_1|$, and as $\rmv{G\bs W_0}{v} \cup W_1 = \rmv{G}{v}$, we obtain $ Ls(G\bs W_0) + |W_0| \ge |\rmv{G}{v}| + Rs(G\bs \rmv{G}{v})$. Finally, as the move $v$ was supposed optimal in $G$, recall that $Ls(G) = |\rmv{G}{v}| + Rs(G\bs \rmv{G}{v})$, yielding to the desired result.
    
\smallskip
    
Now, we prove $(4) :$ Let $v \in W$ be a vertex played by Right in an optimal move in $G\bs W_0$. By definition, we have 

$$Rs(G\bs W_0) =  -|\rmv{G\bs W_0}{v}| + Ls(G\bs W_0 \bs \rmv{G\bs W_0}{v})$$
    
\noindent Denote by $W_1$ the white vertices of $W_0$ taken by playing $v$, i.e. $W_1 = W_0 \cap \rmv{G}{v}$, and by $W_2$ the other white vertices of $W_0$, i.e. $W_2 = W_0 \bs W_1$. 

\noindent We have $ Ls(G\bs W_0 \bs \rmv{G\bs W_0}{v}) = Ls(G\bs W_2 \bs \rmv{G}{v})$.
        
\noindent Now, by induction hypothesis and by applying (1) to $Ls(G\bs W_2 \bs \rmv{G}{v})$, we have 

$$Ls(G\bs W_2 \bs \rmv{G}{v}) \ge Ls(G\bs \rmv{G}{v}) - |W_2|$$

\noindent As $\rmv{G\bs W_0}{v} \subset \rmv{G}{v}$, we have $-|\rmv{G\bs W_0}{v}| \ge - |\rmv{G}{v}|$.

\noindent Putting together the previous inequalities we get $Rs(G\bs W_0) \ge -|\rmv{G}{v}| + Ls(G\bs \rmv{G}{v}) - |W_2|$.

\noindent Since $-|\rmv{G}{v}| + Ls(G\bs \rmv{G}{v})$ is the score obtained by playing $v$ in $G$, we have $Rs(G\bs W_0) \ge Rs(G) - |W_2|$. 

\noindent Finally, as $W_2 \subset W_0$, we obtain $Rs(G\bs W_0) \ge Rs(G) - |W_0|$.
\end{proof}

\begin{exemple}
The graph $G$ of Figure~\ref{fig:exineq} is the sum of a graph and its negative with only two additional black vertices ($v'_7$ and $v'_8$). To evaluate the score of this sum, Left can give these two vertices to Right, and know that on the remaining graph, the score will be $0$ (as being the sum of a graph with its negative). Therefore, by applying $(2)$ and $(3)$ of the above theorem, we deduce that $Ls(G) \ge -2$ and $Rs(G) \ge -2$.

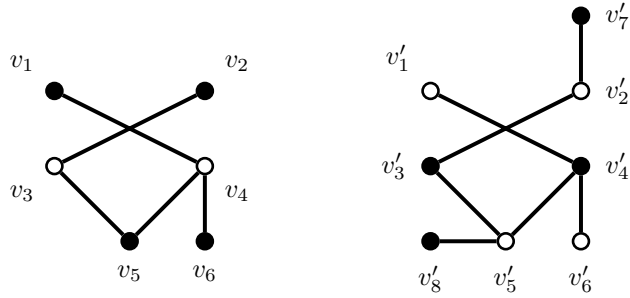
\begin{figure}[ht]
    \centering
\begin{tikzpicture}

  \draw (0,2) node[vertexB](1){} node[above left=0.2] {$v_1$};
  \draw (2,2) node[vertexB](2){} node[above right=0.2] {$v_2$};
  \draw (0,1) node[vertexW](3){} node[below left=0.2] {$v_3$};
  \draw (2,1) node[vertexW](4){} node[below right=0.2] {$v_4$};
  \draw (1,0) node[vertexB](5){} node[below = 0.2] {$v_5$};
  \draw (2,0) node[vertexB](6){} node[below = 0.2] {$v_6$};

  \draw[arete] (1)--(4)--(5)--(3)--(2) (4)--(6);

  \draw (5,2) node[vertexW](1){} node[above left=0.2] {$v'_1$};
  \draw (7,2) node[vertexW](2){} node[right=0.2] {$v'_2$};
  \draw (5,1) node[vertexB](3){} node[left=0.2] {$v'_3$};
  \draw (7,1) node[vertexB](4){} node[right=0.2] {$v'_4$};
  \draw (6,0) node[vertexW](5){} node[below = 0.2] {$v'_5$};
  \draw (7,0) node[vertexW](6){} node[below = 0.2] {$v'_6$};
  \draw (7,3) node[vertexB](7){} node[right = 0.2] {$v'_7$};
  \draw (5,0) node[vertexB](8){} node[below = 0.2] {$v'_8$};

 \draw[arete] (1)--(4)--(5)--(3)--(2)--(7) (4)--(6) (5)--(8);

\end{tikzpicture}

    \caption{Left, by "giving" the two black vertices $v'_7$ and $v'_8$, can ensure a score of 
    $-2$ in this graph since the rest of the graph has a (Left and Right) score of 0.}
    \label{fig:exineq}
\end{figure}
\end{exemple}

\subsection{Twins}

Two vertices are called {\em twins} if they have the same color and the same neighborhood. 
In \bip, twin vertices are necessarily removed together.


\begin{lemma}
\label{lem:twins}
In \bip, if a move removes a vertex $v$ of the graph, then this move also removes all the vertices that are twins of $v$.
\end{lemma}

\begin{proof}
Let $G = (B\cup W,E)$ be a bipartite graph and let $a$ and $a'$ be two twin vertices of $G$. Without loss of generality, suppose that $a$ and $a'$ are white. Consider a move that removes the vertex $a$.

If the move is played from a white vertex $w$, it takes $a$ if and only if it takes all its neighborhood. Thus $N(a) \subseteq N(w)$. But, as $N(a) = N(a')$, we also have $N(a') \subseteq N(w)$. Therefore, $a'$ and $a$ have an empty neighborhood after this move and are won by white.

If the move that removes $a$ is played on a black vertex $b$, then $a$ must be in the neighborhood of $b$. Therefore, since $a$ and $a'$ have the same neighborhood, $a'$ is also in the neighborhood of $b$, which concludes the proof.

\end{proof}










\section{{\sc Pspace}-completeness}\label{sec:pspace}

In this section, we prove that computing the Left score of an instance of \bip\ is {\sc Pspace}-complete.

\

\noindent\decisionpb{\bip}{A bipartite graph $G=(B\cup W, E)$, an integer $k\in \mathbb N$}{Is it true that $Ls(G)\ge k$?}{1}

\begin{theorem}
\bip\ is {\sc Pspace}-complete.
\end{theorem}

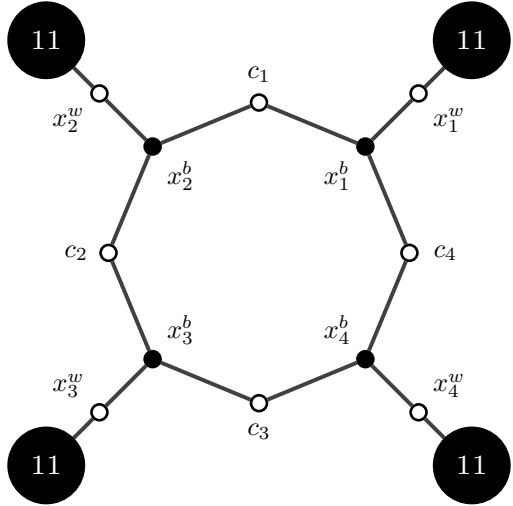
\begin{figure}[h]
    \centering

\begin{tikzpicture}

\draw ( 2.0 , 0.0 ) node[vertexW](s0){} node[right=0.2] {$c_4$};
\draw ( 1.4142135623730951 , 1.4142135623730951 ) node[vertexB](s1){} node[below left=0.1] {$x_1^b$};
\draw ( 1.2246467991473532e-16 , 2.0 ) node[vertexW](s2){}node[above=0.15] {$c_1$} ;
\draw ( -1.414213562373095 , 1.4142135623730951 ) node[vertexB](s3){} node[below right=0.1] {$x_2^b$};
\draw ( -2.0 , 2.4492935982947064e-16 ) node[vertexW](s4){} node[left=0.15] {$c_2$};
\draw ( -1.4142135623730954 , -1.414213562373095 ) node[vertexB](s5){} node[above right=0.1] {$x_3^b$};
\draw ( -3.6739403974420594e-16 , -2.0 ) node[vertexW](s6){} node[below=0.15] {$c_3$};
\draw ( 1.4142135623730947 , -1.4142135623730954 ) node[vertexB](s7){} node[above left=0.1] {$x_4^b$};

\draw ( 2.121320343559643 , 2.121320343559643 ) node[vertexW](e1){} node[below right=0.1] {$x_1^w$} ;
\draw ( -2.1213203435596424 , 2.121320343559643 ) node[vertexW](e3){} node[below left=0.1] {$x_2^w$} ;
\draw ( -2.121320343559643 , -2.1213203435596424 ) node[vertexW](e5){} node[above left=0.1] {$x_3^w$}  ;
\draw ( 2.121320343559642 , -2.121320343559643 ) node[vertexW](e7){} node[above right=0.1] {$x_4^w$} ;

\Vertex[x = 2.8284271247461903,y=2.8284271247461903, color = black, size = 1, fontcolor = white, fontscale = 1.5, label = 11]{k1}
\Vertex[x = -2.82842712474619,y=2.8284271247461903, color = black, size = 1, fontcolor = white, fontscale = 1.5, label = 11]{k3}
\Vertex[x = -2.8284271247461907,y=-2.82842712474619, color = black, size = 1, fontcolor = white, fontscale = 1.5, label = 11]{k5}
\Vertex[x = 2.8284271247461903,y=-2.8284271247461903, color = black, size = 1, fontcolor = white, fontscale = 1.5, label = 11]{k7}

\Edge[](s1)(e1)
\Edge[](s3)(e3)
\Edge[](s5)(e5)
\Edge[](s7)(e7)

\Edge[](k1)(e1)
\Edge[](k3)(e3)
\Edge[](k5)(e5)
\Edge[](k7)(e7)

\Edge[](s0)(s1)
\Edge[](s1)(s2)
\Edge[](s2)(s3)
\Edge[](s3)(s4)
\Edge[](s4)(s5)
\Edge[](s5)(s6)
\Edge[](s6)(s7)
\Edge[](s7)(s0)

\end{tikzpicture}

    \caption{Reduction of $\phi = (X_1 \vee X_2) \wedge (X_2 \vee X_3) \wedge (X_3 \vee X_4) \wedge (X_4 \vee X_1)$ to \bip. Big black vertices represent 11 pendant vertices.}
    \label{fig:reduction}
\end{figure}

\begin{proof}
First, \bip\ is in {\sc Pspace} since a game ends after at most $|B\cup W|$ moves. To prove the completeness, we do a reduction from the game {\sc POS CNF}. An instance of {\sc POS CNF} is a CNF formula $\varphi$ over a set of variables $X=\{X_1,\ldots,X_n\}$ where all the variables are positive. Alice and Bob alternately choose a variable that has not been chosen yet. Variables chosen by Alice are set to true, those chosen by Bob are set to false. Alice wins if and only if $\varphi$ is true at the end.
 The related decision problem is the following one and has been proved to be {\sc Pspace}-complete by Schaefer in 1978 \cite{Schaefer1978}.

\

\noindent\decisionpb{\poscnf}{A set of variable $X$ and a positive CNF formula $\varphi$}{Does Alice win the game \poscnf\ played on $(X, \varphi)$ playing first?}{1}

Let $(X, \varphi)$ be an instance of \poscnf. Without loss of generality, by adding a useless variable, one can assume that $n=|X|$ is even. Denote $\varphi = \underset{j=1}{\overset{m}{\bigwedge}} C_j$ and, for $1 \le j \le m$, denote $C_j = \underset{i=1}{\overset{n_j}{\bigvee}} X^j_i$ where $n_j$ is the number of variables in $C_j$. We construct an instance of \bip\ $G = (B \cup W, E)$ as follows:

\begin{itemize}
    \item For each clause $C_j$, $1 \le j \le m$, we add a white vertex $c_j$.
    \item For each variable $X_i$, $1 \le i \le n$, we add a white and a black vertex ${x_i}^w, {x_i}^b$, and $m+2n-1$ black vertices $v^i_1, \dots v^i_{m+2n-1}$. 
    \item We finally add all the edges $({x_i}^w, {x_i}^b)$ and $({x_i}^w,v^i_k)$ for $1 \le i \le n$ and $1 \le k \le m+2n$; and the edges $({x_i}^b,c_j)$ whenever $X_i \in C_j$.
\end{itemize}

See Figure \ref{fig:reduction} for an illustration of the construction for the formula $\varphi = (X_1 \vee X_2) \wedge (X_2 \vee X_3) \wedge (X_3 \vee X_4) \wedge (X_4 \vee X_1)$. The big black vertices represent the $m+2n-1=11$ black leaves attached to each vertex ${x_i}^w$.

Before proving the reduction, we prove that there are optimal strategies for Left or Right in $G$ consisting in playing only vertices ${x_i}^w$ or ${x_i}^b$ (for $1 \le i \le n$). 
For $1\leq i \leq n$, let $B_i$ be the set of $m+2n-1$ black vertices $v^i_k$. We call such set a {\em bag}. All these vertices are twins and by Lemma \ref{lem:twins} will be removed at the same time, together with the vertex $x_i^w$, scoring at least $2n+m$ points. However, a move can remove at most one bag at a time, and actually any move except the moves on vertices $c_j$ removes a bag.
As a consequence, if at some point Right chooses to play on a vertex $c_j$ whereas a bag is still present, Left will manage to take $n/2+1$ bags and will remove at least
$\left ( \frac{n}{2} + 1 \right ) (2n+m) = n^2 + \frac{nm}{2} + 2n + m$ vertices.
On the contrary, if Right always takes a bag when one is left, Left can take at most $n/2$ bags, plus maybe the $n+m$ remaining vertices, removing at most $\left ( \frac{n}{2} \right ) (2n+m)+n+m = n^2 + \frac{nm}{2} +n +m$ vertices.

Thus, until all the bags are removed, Right must always choose to play on a vertex $x_i^w$. When the bag $B_i$ is removed (by Left or Right), the only other vertices that can be removed in the same move are $x_i^w$, $x_i^b$ and maybe some clause vertices, but no vertices of other variables can be removed. Thus, if Left would have play on a vertex $v^i_k$ in a bag, the vertex $x_i^b$ is still present in the graph. Since $\rmv{G,v^i_k} \subset \rmv{G,{x_i}^b}$, by Theorem \ref{thm:includedmoves}, playing on $x_i^b$ is always better.
Finally, we can suppose that both Left and Right plays on vertices $x_i^b$ and $x_i^r$ until all the bags are removed. But then the graph is empty since only clauses variables could remain and they form an independent set.

Note that the final Left score is exactly the difference between the number of vertices $c_j$ removed by Left (by taking one adjacent vertex $x_i^b$ of $c_j$) and thus isolated by Right (by removing all the black vertices $x_i^b$ adjacent to $c_j$).

Assume Alice has a winning strategy in \poscnf\ on $(X,\varphi)$. Then Left follows the same strategy than Alice: when Alice chooses variable $X_i$, Left removes vertex $x_i^b$. If Right answers vertex $x_{i'}^w$, we assume that Bob plays variable $X_{i'}$ and so on. At the end, Left will have take one neighbour of each clause vertex and thus will have removed the $m$ clause vertices. At the end, $Ls(G)\geq m$.

Similarly, if Left has a strategy to have $Ls(G)\geq m$, it means that she can choose vertices $x_i^b$ in such a way that she will touch all the clause variables. If Alice follows the same strategy, $\varphi$ will be true at the end of the game.

\end{proof}



\begin{remark}

 In \bip, the question was to know whether $Ls(G) \ge k$. Since free points can be given to a player by adding isolated vertices in the graph, the problem is still {\sc Pspace}-complete with $k = 0$ by adding $k$ white vertices if $k>0$ or $-k$ black vertices if $k<0$. Thus, determining which player has a winning strategy is also {\sc  Pspace}-complete.
\end{remark}

\begin{corollaire}
{\sc influence} is {\sc Pspace}-complete, even for bipartite graphs.
\end{corollaire}

\section{A sufficient condition to have a draw}
\label{sec:symetrie}

Identifying {\em draw} games, i.e. games that are equivalent to $0$, is essential in the study of a game in Milnor's universe since they can be removed from any sum. However, as seen in the previous section, determining in general if a game is a draw is hard (since it is already hard to determine if $Ls(G)=0$). In this section, we give a sufficient condition to have a draw when the graph has some symmetries and we apply it to hypercubes, cylinders and torus.

\subsection{Sufficient condition}

An {\em automorphism} $\varphi$ of a graph $G$ is a function from $V(G)$ to $V(G)$ that preserves edges, i.e. if $(u,v) \in E(G)$ if and only if $(\varphi(u),\varphi(v)) \in E(G)$.
It is said to be involutive if $\varphi(\varphi(u)) = u$ for all $u \in V(G)$.
We denote the usual distance (i.e the length of the shortest path) between two vertices $u$ and $v$ by $d(u,v)$. If $u$ and $v$ are not in the same connected component, we let $d(u,v)=\infty$.

\begin{definition}
Let $G = (B\cup W,E)$ be a bipartite graph. An automorphism $\varphi$ of $G$ is a {\it BW-automorphism} if it is involutive, exchanging the colours of the vertices, and if for any $v \in V(G)$, $d(v, \varphi(v)) \ge 3$.
\end{definition}

\begin{theorem}\label{thmsym}
Let $G = (B\cup W,E)$ be a bipartite graph that has a BW-automorphism, then $G=0$.
\end{theorem}

\begin{proof}
 
Let $G = (B\cup W,E)$ a bipartite graph, and let $\varphi$ be a BW-automorphism of $G$. To prove that $G=0$, we just need to prove that $LS(G) = RS(G) = 0$. By symmetry of $G$ and $\varphi(G)$, we actually just need to prove that $Ls(G)\leq 0$, that is to prove that the second player can achieve to take as many vertices as the second player. The strategy of the second player will consist in always playing the image by $\varphi$ of the vertex played by the first player.

To prove that this strategy works, it is sufficient to prove that $\rmv{G}{\varphi(u)} = \varphi(\rmv{G}{u})$ for all $u \in G$ and that $\rmv{G}{\varphi(u)} \cap \rmv{G}{u} = \emptyset$. Indeed, if these two conditions are satisfied, then the second player can always remove the image by $\varphi$ of the moved played by the first player and will take the same number of vertices. Then after two moves, the graph will still have a BW-automorphism and the second player can repeat his strategy.

Let $u \in V$ and $v \in \rmv{G}{\varphi(u)}$. If $v$ is a neighbour of $\varphi(u)$, we have :

\begin{align*}
    (v,\varphi(u)) \in E
  & \Leftrightarrow (\varphi(v), \varphi ( \varphi(u))) \in E \mbox{ as $\varphi$ respects edges}\\
  & \Leftrightarrow (\varphi(v), u) \in E \mbox{ as $\varphi$ is involutive.}
 \end{align*}

 If $v$ is not a neighbour of $\varphi(u)$, $v$ is isolated by playing $\varphi(u)$, i.e. each neighbour $w$ of $v$ is also a neighbour of $\varphi(u)$. As $\varphi$ preserves edges, the neighbours of $\varphi(v)$ are the images by $\varphi$ of the neighbours of $v$. Therefore, they are all removed by playing $\varphi(\varphi(u)) = u$. This proves that $\varphi(v)$ is isolated by playing $u$, and thus is in $\rmv{G}{u}$.
 
 In both cases, $\varphi(v) \in \rmv{G}{u}$ and thus $v=\varphi(\varphi(v)$ is in $\varphi(\rmv{G}{u})$. This proves $\rmv{G}{\varphi(u)} \subseteq \varphi(\rmv{G}{u})$ for all $u \in V$. 
 
 To prove the other direction, we have the following equivalences:
 \begin{align*}
    &\forall u \in V, \varphi( \rmv{G}{\varphi(u)} ) \subseteq \varphi( \varphi(\rmv{G}{u}) )  \\
    \Leftrightarrow  &\forall u \in V, \varphi( \rmv{G}{\varphi(u)} ) \subseteq \rmv{G}{u} \mbox{ as $\varphi$ is involutive} \\
    \Leftrightarrow &\forall v \in V, \varphi( \rmv{G}{\varphi(\varphi(v))}) \subseteq \rmv{G}{\varphi(v)} \mbox{ as $\varphi$ is a bijection, by taking $v = \varphi(u)$ }\\
    \Leftrightarrow &\forall v \in V, \varphi( \rmv{G}{v} ) \subseteq \rmv{G}{\varphi(v)} \mbox{ as $\varphi$ is involutive.}
 \end{align*}
 
\noindent Finally, we have $\rmv{G}{\varphi(u)} = \varphi(\rmv{G}{u})$ for all $u \in V$. 

A move in \bip\ removes vertices at distance at most 2, and the vertices at distance 2 are removed only if they are made isolated. Assume by contradiction there exists $v \in \rmv{G}{\varphi(u)} \cap \rmv{G}{u}$. Then $d(v,u)\le 2$ and $d(v,\varphi(u))\leq 2$. Since $d(u,\varphi(u))\geq 3$, we actually have $d(u,\varphi(u))=3$ (since $u$ and $\varphi(u)$ does not have the same colour, and without loss of generality, we can assume that $d(v,u)=2$ and $d(v,\varphi(u))=1$. But then $v$ cannot be isolated when playing $u$ since it will still be connected to $\varphi(u)$. Thus $\rmv{G}{\varphi(u)} \cap \rmv{G}{u}=\emptyset$.

\end{proof}

\subsection{Applications}

In this subsection, we apply Theorem \ref{thmsym} to hypercubes, torus and cylinders.

\begin{figure}
    \centering
\scalebox{0.6}{\begin{tikzpicture}

 \Vertex[x=0,y=0, color = white, size = 1, fontcolor = black, fontscale = 1.6, label = $0000$]{0000}
 
 \Vertex[x=3,y=0, color = black, size = 1, fontcolor = white, fontscale = 1.6, label = $0001$]{0001}

 \Vertex[x=3,y=3, color = white, size = 1, fontcolor = black, fontscale = 1.6, label = $0011$]{0011}
 
 \Vertex[x=0,y=3, color = black, size = 1, fontcolor = white, fontscale = 1.6, label = $0010$]{0010}

  \Vertex[x=1.5,y=1, color = black, size = 1, fontcolor = white, fontscale = 1.6, label = $0100$]{0100}
 
 \Vertex[x=4.5,y=1, color = white, size = 1, fontcolor = black, fontscale = 1.6, label = $0101$]{0101}

 \Vertex[x=4.5,y=4, color = black, size = 1, fontcolor = white, fontscale = 1.6, label = $0111$]{0111}
 
 \Vertex[x=1.5,y=4, color = white, size = 1, fontcolor = black, fontscale = 1.6, label = $0110$]{0110}

     \Edge[bend = 0](0000)(0001)
     \Edge[bend = 0](0010)(0011)
     \Edge[bend = 0](0000)(0010)
     \Edge[bend = 0](0011)(0001)
     
     \Edge[bend = 0, style = dashed](0100)(0101)
     \Edge[bend = 0](0110)(0111)
     \Edge[bend = 0, style = dashed](0100)(0110)
     \Edge[bend = 0](0111)(0101)
     
     \Edge[bend = 0, style = dashed](0000)(0100)
     \Edge[bend = 0](0010)(0110)
     \Edge[bend = 0](0001)(0101)
     \Edge[bend = 0](0011)(0111)

 \Vertex[x=7,y=0, color = black, size = 1, fontcolor = white, fontscale = 1.6, label = $1000$]{1000}
 
 \Vertex[x=10,y=0, color = white, size = 1, fontcolor = black, fontscale = 1.6, label = $1001$]{1001}

 \Vertex[x=10,y=3, color = black, size = 1, fontcolor = white, fontscale = 1.6, label = $1011$]{1011}
 
 \Vertex[x=7,y=3, color = white, size = 1, fontcolor = black, fontscale = 1.6, label = $1010$]{1010}

  \Vertex[x=8.5,y=1, color = white, size = 1, fontcolor = black, fontscale = 1.6, label = $1100$]{1100}
 
 \Vertex[x=11.5,y=1, color = black, size = 1, fontcolor = white, fontscale = 1.6, label = $1101$]{1101}

 \Vertex[x=11.5,y=4, color = white, size = 1, fontcolor = black, fontscale = 1.6, label = $1111$]{1111}
 
 \Vertex[x=8.5,y=4, color = black, size = 1, fontcolor = white, fontscale = 1.6, label = $1110$]{1110}

     \Edge[bend = 0](1000)(1001)
     \Edge[bend = 0](1010)(1011)
     \Edge[bend = 0](1000)(1010)
     \Edge[bend = 0](1011)(1001)
     
     \Edge[bend = 0, style = dashed](1100)(1101)
     \Edge[bend = 0](1110)(1111)
     \Edge[bend = 0, style = dashed](1100)(1110)
     \Edge[bend = 0](1111)(1101)
     
     \Edge[bend = 0, style = dashed](1000)(1100)
     \Edge[bend = 0](1010)(1110)
     \Edge[bend = 0](1001)(1101)
     \Edge[bend = 0](1011)(1111)
     
     \Edge[bend = 270, color = black](0000)(1000)
     \Edge[bend = -270, color = black](0010)(1010)
     \Edge[bend = 270, color = black](0001)(1001)
     \Edge[bend = -270, color = black](0011)(1011)
     
     \Edge[bend = 270, color = black](0100)(1100)
     \Edge[bend = -270, color = black](0110)(1110)
     \Edge[bend = 270, color = black](0101)(1101)
     \Edge[bend = -270, color = black](0111)(1111)

\end{tikzpicture}}
    \caption{Hypercube $H_4$}
    \label{fig:hc}
\end{figure}
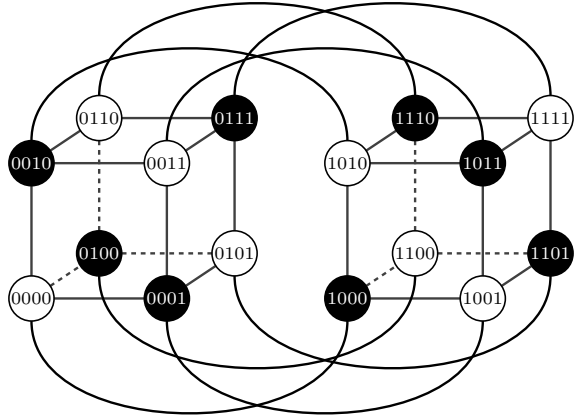

\begin{definition}[Hypercube]
We denote by $H_n$ be the n-dimensional hypercube. The set of vertices of $H_n$ is $\{0,1\}^n$ and two vertices are adjacent if and only if they differ on exactly one digit. We set black vertices to be the vertices with an odd number of nonzero digits. See Figure \ref{fig:hc} for a representation of $H_4$.
\end{definition}

\begin{proposition}
For $n\ge 3$, $H_n=0$.
\end{proposition}

\begin{proof}
By Theorem~\ref{thmsym}, we just need to find a BW-automorphism.
Let $n\geq 3$ and 

$\varphi : 
\begin{array}{ccc}
     \{0,1\}^n & \to & \{0,1\}^n \\
     (\eps_1, \dots, \eps_n) & \mapsto & (1-\eps_1, 1- \eps_2, 1-\eps_3, \eps_4, \dots, \eps_n)
\end{array}$.

The automorphism $\varphi$ is involutive and  $d(u,\varphi(u))=3$ for any vertex $u$, which implies that $\varphi$ exchanges black and white vertices and. Thus $\varphi$ is a BW-automorphism and $H_n=0$.
\end{proof}

\begin{definition}[Cylinder]
Let $n,m$ be two positive integers with $n$ even. We denote $C_{n,m}$ the cylinder defined as follows. The $n\times m$ vertices of $C_{n,m}$ are denoted by $v_{i,j}$ with $1\leq i \leq n$ and $1\leq j\leq m$. Vertices $v_{i,j}$ and $v_{i',j'}$ are adjacent if and only if $i=i'$ and $j=j' \pm 1$ or $j=j'$ and $i=i'\pm 1 \bmod n$. We set black vertices to be the vertices with $i+j=0\bmod 2$.
\end{definition}

\begin{definition}[Torus]
Let $n,m$ two even positive integers. We denote by $T_{n,m}$ the torus defined as follows. 
The $n\times m$ vertices of $T_{n,m}$ are denoted by $v_{i,j}$ with $1\leq i \leq n$ and $1\leq j\leq m$. Vertices $v_{i,j}$ and $v_{i',j'}$ are adjacent if and only if $i=i'$ and $j=j' \pm 1 \bmod m$ or $j=j'$ and $i=i'\pm 1 \bmod n$. We set black vertices to be the vertices with $i+j=0\bmod 2$.

\end{definition}

\begin{figure}[h]
    \centering
\scalebox{0.8}{\begin{tikzpicture}

  \Vertex[x=0,y=-0, color = black, size = 1, fontcolor = white, fontscale = 1.6, label = $v_{1,1}$]{11}

 \Vertex[x=2,y=-0, color = white, size = 1, fontcolor = black, fontscale = 1.6, label = $v_{1,2}$]{12}

  \Vertex[x=4,y=-0, color = black, size = 1, fontcolor = white, fontscale = 1.6, label = $v_{1,3}$]{13}

 \Vertex[x=6,y=-0, color = white, size = 1, fontcolor = black, fontscale = 1.6, label = $v_{1,4}$]{14}

  \Vertex[x=8,y=-0, color = black, size = 1, fontcolor = white, fontscale = 1.6, label = $v_{1,5}$]{15}

 \Vertex[x=10,y=-0, color = white, size = 1, fontcolor = black, fontscale = 1.6, label = $v_{1,6}$]{16}

 \Vertex[x=0,y=-2, color = white, size = 1, fontcolor = black, fontscale = 1.6, label = $v_{2,1}$]{21}

  \Vertex[x=2,y=-2, color = black, size = 1, fontcolor = white, fontscale = 1.6, label = $v_{2,2}$]{22}

 \Vertex[x=4,y=-2, color = white, size = 1, fontcolor = black, fontscale = 1.6, label = $v_{2,3}$]{23}

  \Vertex[x=6,y=-2, color = black, size = 1, fontcolor = white, fontscale = 1.6, label = $v_{2,4}$]{24}

 \Vertex[x=8,y=-2, color = white, size = 1, fontcolor = black, fontscale = 1.6, label = $v_{2,5}$]{25}

  \Vertex[x=10,y=-2, color = black, size = 1, fontcolor = white, fontscale = 1.6, label = $v_{2,6}$]{26}

  \Vertex[x=0,y=-4, color = black, size = 1, fontcolor = white, fontscale = 1.6, label = $v_{3,1}$]{31}

 \Vertex[x=2,y=-4, color = white, size = 1, fontcolor = black, fontscale = 1.6, label = $v_{3,2}$]{32}

  \Vertex[x=4,y=-4, color = black, size = 1, fontcolor = white, fontscale = 1.6, label = $v_{3,3}$]{33}

 \Vertex[x=6,y=-4, color = white, size = 1, fontcolor = black, fontscale = 1.6, label = $v_{3,4}$]{34}

  \Vertex[x=8,y=-4, color = black, size = 1, fontcolor = white, fontscale = 1.6, label = $v_{3,5}$]{35}

 \Vertex[x=10,y=-4, color = white, size = 1, fontcolor = black, fontscale = 1.6, label = $v_{3,6}$]{36}

 \Vertex[x=0,y=-6, color = white, size = 1, fontcolor = black, fontscale = 1.6, label = $v_{4,1}$]{41}

  \Vertex[x=2,y=-6, color = black, size = 1, fontcolor = white, fontscale = 1.6, label = $v_{4,2}$]{42}

 \Vertex[x=4,y=-6, color = white, size = 1, fontcolor = black, fontscale = 1.6, label = $v_{4,3}$]{43}

  \Vertex[x=6,y=-6, color = black, size = 1, fontcolor = white, fontscale = 1.6, label = $v_{4,4}$]{44}

 \Vertex[x=8,y=-6, color = white, size = 1, fontcolor = black, fontscale = 1.6, label = $v_{4,5}$]{45}

  \Vertex[x=10,y=-6, color = black, size = 1, fontcolor = white, fontscale = 1.6, label = $v_{4,6}$]{46}

     \Edge[bend = 0](11)(12)
     \Edge[bend = 0](12)(13)
     \Edge[bend = 0](13)(14)
     \Edge[bend = 0](14)(15)
     \Edge[bend = 0](15)(16)
     \Edge[bend = -200, style = dashed](16)(11)
     
     \Edge[bend = 0](21)(22)
     \Edge[bend = 0](22)(23)
     \Edge[bend = 0](23)(24)
     \Edge[bend = 0](24)(25)
     \Edge[bend = 0](25)(26)
     \Edge[bend = -200, style = dashed](26)(21)

     \Edge[bend = 0](31)(32)
     \Edge[bend = 0](32)(33)
     \Edge[bend = 0](33)(34)
     \Edge[bend = 0](34)(35)
     \Edge[bend = 0](35)(36)
     \Edge[bend = -200, style = dashed, style = dashed](36)(31)
     
     \Edge[bend = 0](41)(42)
     \Edge[bend = 0](42)(43)
     \Edge[bend = 0](43)(44)
     \Edge[bend = 0](44)(45)
     \Edge[bend = 0](45)(46)
     \Edge[bend = -200, style = dashed](46)(41)     

     \Edge[bend = 0](11)(21)
     \Edge[bend = 0](21)(31)
     \Edge[bend = 0](31)(41)
     \Edge[bend = -200](41)(11)
     
     \Edge[bend = 0](12)(22)
     \Edge[bend = 0](22)(32)
     \Edge[bend = 0](32)(42)
     \Edge[bend = -200](42)(12)
     
     \Edge[bend = 0](13)(23)
     \Edge[bend = 0](23)(33)
     \Edge[bend = 0](33)(43)
     \Edge[bend = -200](43)(13)
     
     \Edge[bend = 0](14)(24)
     \Edge[bend = 0](24)(34)
     \Edge[bend = 0](34)(44)
     \Edge[bend = -200](44)(14)
     
     \Edge[bend = 0](15)(25)
     \Edge[bend = 0](25)(35)
     \Edge[bend = 0](35)(45)
     \Edge[bend = -200](45)(15)

      \Edge[bend = 0](16)(26)
     \Edge[bend = 0](26)(36)
     \Edge[bend = 0](36)(46)
     \Edge[bend = -200](46)(16)

\end{tikzpicture}}
    \caption{Cylinder $C_{4,6}$ (Torus $T_{4,6}$ with the dashed edges)}
    \label{fig:torus}
\end{figure}
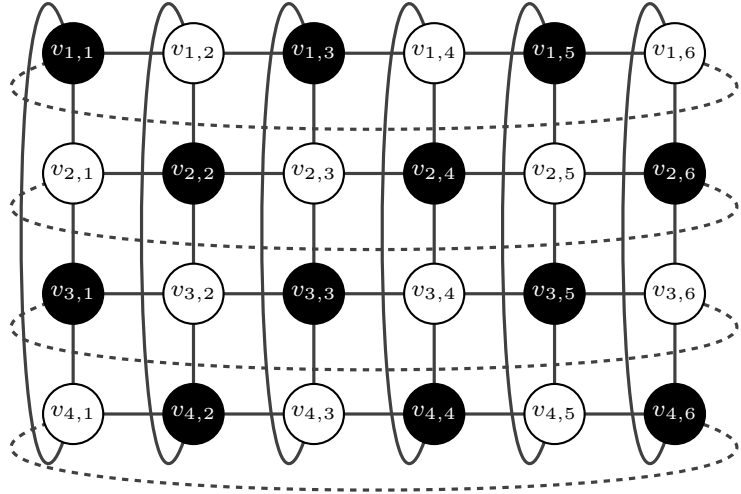

See Figure \ref{fig:torus} for a representation of $C_{4,6}$ and $T_{4,6}$. Note that a torus is a cylinder where the two extremal columns are connected.

\begin{proposition}
If $n \ge 4$ or $m\ge 4$, then $T_{n,m}=0$.

If $n \ge 4$ and $m$ is even or if $n = 4k+2$ for a positive integer k, we have $C_{n,m}= 0$.
\end{proposition}

\begin{proof}
By symmetry of the torus, we can assume that $n\geq 4$.

Assume first that $n=2\bmod 4$. Consider the following morphism on the vertices of $C_{n,m}$ or $T_{n,m}$: $$\varphi:  v_{i,j} \mapsto  v_{((i+\frac{n}{2}) \mod n),j)}.$$
It is sending vertices to the symmetric vertex of that is in the same cycle. Thus $\varphi$ is involutive. Since $n=2\bmod 4$, the colour of $u$ and $\varphi(u)$ are opposite. Since $n\geq 6$, they are distance at least 3. Thus $\varphi$ is a BW-automorphism and $C_{n,m}=T_{n,m}=0$.

Assume now that $n=0\bmod 4$ and that both $n$ and $m$ are even.
Consider now the morphism $$\psi :
     v_{i,j}  \mapsto  v_{((i+\frac{n}{2}) \mod n),(m+1-j)}.$$
One can check that $\psi$ is an automorphism.
Furthermore, 
$$\psi(\psi(v_{i,j})) = \psi(v_{((i+\frac{n}{2}) \mod n),(m+1-j)}) = v_{((i+\frac{n}{2} + \frac{n}{2}) \mod n),(m+1-(m+1-j))} = v_{i,j}$$ and thus $\psi$ is involutive.

Since $i+j$ and $i+\frac{n}{2}+m+1-j$ have not the same parity, $\psi$ is exchanging black and white vertices.

Finally,  $d(v_{i,j}, \psi(v_{i,j})) \ge \frac{n}{2} + 1\geq 3$. Indeed, to go from $v_{i,j}$ to $\psi(v_{i,j})$ one needs already $n/2$ steps to go to the same line, and at least one more step to go to the same column since $m+1-j\neq j$ ($m$ is even).

Thus $\psi$ is a BW-automorphism and
by Theorem \ref{thmsym} we can conclude that $T_{n,m}=C_{n,m} = 0$.
\end{proof}

\subsection{Limits of Theorem \ref{thmsym}}
Despite Theorem \ref{thmsym} gives a good condition to find some draw graphs, we prove now that deciding if a BW-automorphism exists in a graph is not solvable in polynomial time, unless the graph isomorphism problem is.

\begin{definition}[Graph Isomorphism problem]
We recall here that given two graphs $G_1$ and $G_2$, determining if $G_1$ and $G_2$ are isomorphic is not known to be solvable in polynomial time. Problems that are known to be at least as difficult as this problems are called {\sc GI}-Hard problems.
\end{definition}

Note that the graph isomorphism problem is known to be computable in quasi-polynomial time \cite{lasz}. 

\begin{theorem}
Let $G = (B\cup W,E)$ be a bipartite graph, deciding if $G$ admits a BW-automorphism is {\sc GI}-Hard.
\end{theorem}

\begin{proof}

We give a reduction from the graph isomorphism problem. Let $G_1 = (V_1,E_2), G_2 = (V_2,E_2)$ be two graphs. 

We construct an instance $H=(B\cup W,E)$ of \bip\ as follows. The graph $H$ will be divided into two bipartite graphs, $H_1$ and $H_2$. Let $i\in\{1,2\}$. The graph $H_i$ is the incidence graph of $G_i$ where a leaf is added to each vertex of $H_i$ that corresponds to a vertex of $G_i$.
More precisely, the vertex set of $H_i$ is $V_i\cup E_i\cup V'_i$ where $V'_i$ is a copy of $V_i$. Two vertices $u$ and $v$ are adjacent in $H_i$ if and only if one of the following items is satisfied
\begin{itemize}
    \item $u \in V_i$ and $v$ is the copy of $u$ in $V'_i$, or
    \item $u \in V_i$, $v\in E_i$ and $v$ is an edge incident to $u$ in $G_i$.
\end{itemize}

The graph $H_i$ is bipartite with parts $V_i$ and $E_i\cup V'_i$. Note that the only vertices of degree $1$ are the vertices in $V'_i$.
Then $H=H_1\cup H_2$. We choose as set of black vertices the set $B=V_1\cup E_2\cup V'_2$ and as set of white vertices the set $W=V_2 \cup E_1 \cup V'_1$.

We now prove that $G_1$ and $G_2$ are isomorphic if and only if there exists a BW-automorphism of $H$.

First suppose that $G_1$ and $G_2$ are isomorphic. Let $\phi$ be an isomorphism from $G_1$ to $G_2$. We construct a BW-automorphism $\varphi$ of $H$ as follows. Vertices of $V_1$ and $V_2$ are exchanged as in $\phi$. Since $\phi$ preserves edges, a vertex of $H_1$ corresponding to an edge $(u,v)$ of $G_1$ is send to the vertex of $H_2$ corresponding to the edge $(\phi(u),\phi(v))$ of $G_2$, and conversely for the vertices of $H_2$ corresponding to the edges of $G_2$.
Finally, a vertex of $V'_1$ (respectively $V'_2$) that is a copy of a vertex $u$ of $V_1$ (resp. $V_2$) is send to the copy of $\phi(u)$ in $V'_2$ (resp. to the copy of $\phi^{-1}(u)$ in $V'_1$). Clearly, $\varphi$ is an automorphism that is involutive, is inversing the colour of the vertices, and since it is sending a vertex to a vertex not in the same component, the distance between $u$ and $\varphi(u)$ is infinite. Thus, $\varphi$ is a BW-automorphism of $H$.

Assume now that a BW-automorphism $\varphi$ of $H$ exists. 
Since $V'_1$ and $V'_2$ are the only vertices of degree 1 of $H_1$ and $H_2$ respectively, they must be stable by $\varphi$. To respect the exchange of colours, $\varphi$ must induced a bijection between $V'_1$ and $V'_2$. Indeed, the set $V_1$ corresponds to the unique neighbours of the vertices of $V'_1$ and thus must be sent to the unique neighbours of the vertices of $V'_2$ that correspond to the set $V_2$. Finally, the vertices of $E_1$ can only be sent to the vertices of $E'_2$. 

Let $\phi:V_1\to V_2$ be the restriction of $\varphi$ to $V_1$. Then $\phi$ is an isomorphism of $G_1$. Indeed, $\varphi$ is involutive so $\phi$ is a bijection. Consider $u$ and $v$ to neighbours in $G_1$ and $e=(u,v)$. Then $\varphi(e)$ must be incident to $\phi(u)$ and $\phi(v)$ and thus $\phi(u)$ and $\phi(v)$ are adjacent in $G_2$. The reverse can be proved in the same way.

To conclude, $G$ admits a BW-automorphism if and only if $G_1$ and $G_2$ are isomorphic, which proves that determining whether a bipartite graph admits a BW-automorphism or not is {\sc GI}-Hard.

\end{proof}

\section{Segments}\label{sec:segments}

\subsection{State of the art}

Recall that segments, as defined in~\cite{duchene}, correspond to instances of \bip\ played on paths. Despite their simple graph structure, their resolution is not obvious.

\begin{definition}[Duchêne et al \cite{duchene}]

The class $C^{segment}$ represents the class of all \bip\ graphs that are paths. 

\end{definition}

When playing on a segment, the graph may be split into two smaller segments. This led the study to consider the class of disjunctive sum of segments (also corresponding to unions of paths), in order to play on a family that is closed by a move operation. 

\begin{definition}

The class $C$ represents the class of all bipartite graphs that are finite disjunctive sums of elements in $C^{segment}$.

\end{definition}

\begin{figure}[ht]
    \centering
    
\begin{tikzpicture}
  
  \node[vertexW](0) at (-0.5,0.5) {};
  \node[vertexB](1) at (0,0) {};
  \node[vertexW](2) at (0.5,0.5) {};
  \node[vertexB](3) at (1,0) {};
  \node[vertexW](4) at (1.5,0.5) {};
  \node[vertexB](5) at (2,0) {};
  \node[vertexW](6) at (2.5,0.5) {};

  \draw[arete] (0)--(1)--(2)--(3)--(4)--(5)--(6);

\end{tikzpicture}
    
    \caption{Segment $S_{-7}$}
    \label{segment}
\end{figure}
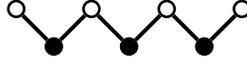

A segment of $n$ vertices starting with a black (resp. white) vertex will be denoted by $S_n$ (resp. $S_{-n}$). See Figure~\ref{segment} for an example.

On the class $C$ without isolated vertex (as such vertices can be immediately removed from any instance), remark that a move removes between $2$ and $5$ vertices. A move that removes $k$ vertices will be called a $k$-move. $4$-moves are the moves made by playing the vertex at distance $2$ from the extremity of a segment. $5$-moves can only be made by Right on $S_{-5}$ and by Left on $S_{5}$, by playing the middle vertex.

By definition of the negative of a game, remark also that $S_{-n} = -S_n$ for all $n \in \N$. Moreover, by symmetry, we have $S_{2n} = S_{-2n}$.\\

On the class $C^{segment}$, Duchêne et al. \cite{duchene} showed that the first player always wins. More precisely, the scores are bounded according to the following result:

\begin{theorem}[Duchêne et al \cite{duchene}]

\label{thm6}

Let $S_n \in C^{segment}$. We have:

\begin{align*}
  -5 \le Rs(S_n) < 0 < Ls(S_n) \le 5
\end{align*}

\end{theorem}

The main idea of the proof of this theorem is that, by playing in the middle vertex of a segment, the first player can create two segments having a number of vertices that is very close. Consequently, this player can ensure to lose at most two vertices in the remaining game. As his first move scores three vertices, it ensures a win for the first player.\\

In addition to this result, there are certain values of $n$ for which the exact value of $Ls(S_n)$ is known (e.g. when $n\equiv 1 \bmod 4$). Yet, for the majority of them, the exact value of the score is an open problem. In the study \cite{duchene}, the first 80 values of $Ls$ and $Rs$ are computed on segments, yielding Figure~\ref{table38} (depicting the first 38 values). The structure of these results led to the following conjecture:

\begin{conjecture}[Duchêne et al \cite{duchene}]\label{conj:periodicity}
The sequences $Ls(S_n)$ and $Rs(S_n)$ are ultimately periodic.
\end{conjecture}

\begin{figure}
    \centering
\scalebox{0.8}{
\begin{tabular}{|c||r|r|r|r|r|r|r|r|r|r|r|r|r|r|r|r|r|r|r|r|}
\hline
n &  1 & 2& 3& 4& 5& 6& 7& 8& 9& 10& 11& 12& 13&  14&   15 &  16 &   17 & 18 & 19 \\
\hline
Ls($S_n$)&1&2 & 3 & 4 & 5 & 2 & 1 & 2 & 3 & 2 & 1 & 2 & 3 &  4 &  3 &  2 & 3 & 2 & 3\\
Rs($S_n$)&1&-2 & -3 & -4 & -1 & -2 & -3 & -2 & -1 & -2 & -3 & -2 & -1 & -4 &  -3 &  -2&  -1 & -2 & -3\\
\hline
\hline
n &  20&  21&  22&  23&  24&   25&26& 27& 28& 29& 30& 31& 32& 33& 34 & 35 & 36 & 37 & 38\\
\hline
Ls($S_n$) &  2 &  5 &   4 &   3 &   2 &  3 & 2 & 3 & 2 & 5 & 4 & 3 & 2 & 3 & 2 & 3 & 2 & 5 & 2 \\
Rs($S_n$) &  -2 &  -1  &   -4 &   -3 &   -2 &   -1 & -2 & -3 & -2 & -1 & -4 & -3 & -2 & -1 & -2 & -3 & -2 & -1 & -2\\
\hline
\end{tabular}

} 

\caption{Left and Right score for segments in \cite{duchene} }
    \label{table38}
\end{figure}

Such a periodicity of the score would directly induce a polynomial time algorithm to solve the game on segments. \\

Concerning sums of segments, Theorem~\ref{thm6} does not hold when extended to the class $C$. The best bounds that are currently known are given by the result below.

\begin{proposition}[Duchêne et al \cite{duchene}]

Let $G \in C$. We have:
\begin{itemize}
    \item $-4 \le Rs(G) \le 0 \le Ls(G) \le 4$ if all the segments of $G$ have even size,  
    \item $-5 \le Rs(G) \le 1$ and $-1 \le Ls(G) \le 5$ if all the segments of $G$ have even size except one. 
\end{itemize}

\label{prop:allevent}
\end{proposition}

In addition, the bounds in this result are tight. For example, consider the sum $G=S_9+S_2$ for which $Rs(G)=1$. This case is the smallest non-trivial graph in $C$ where the first player does not win.\\

\subsection{Mean and temperature of segments}

As mentioned above, sums of segments have not been considered when there are at least two odd segments in the sum. A way to address this issue is to consider the notion of mean and temperature of segments, by applying Corollary\ref{coro:sumtemp}.  \\

We first deal with the computation of the mean of segments. Recall first that $m_L(S_n)=m_R(S_n)$ by Theorem~\ref{thmmoy}. Consequently, we will denote $m_n = m_L(S_n) = m_R(S_n)$.

\begin{theorem}
\label{thm moy segment}
Let $n \in \Z$.
We have $ \left\{
    \begin{array}{lll}
        m_n = 0 & \mbox{if } n \equiv 0 \mod 2 \\
        0 \le m_n \le 1 & \mbox{if } n \equiv 1 \mod 2 \mbox{ and $n \ge 0$}  \\
        -1 \le m_n \le 0 & \mbox{if } n \equiv 1 \mod 2 \mbox{ and $n \le 0$} \\
    \end{array}
\right.$
\end{theorem}

\begin{proof}

First, suppose that n is even. We know that $S_{n} = S_{-n}$. Therefore, if $p$ is even, $\frac{Ls(pS_{n})}{p} = \frac{Ls(0)}{p} = 0$ and if $p$ is odd, we have $\frac{Ls(pS_{n})}{p} = \frac{Ls(S_{n})}{p}$. As $ - \frac{5}{p} \le \frac{Ls(S_{n})}{p} \le \frac{5}{p}$, we have $\frac{Ls(S_{n})}{p} \underset{p \to \infty}{\to} 0$.

\noindent Finally, $\underset{p \to \infty}{\lim} \frac{Ls(pS_{n})}{p} = 0$ and $m_{n} = 0$.\\

Now suppose that $n$ is odd. Since $S_{-n} = -S_n $, and as the mean is linear, we can assume, without loss of generality, that $n \ge 0$. We will show that $0 \le m_n \le 1$.

\noindent If $n=1$, we have $Ls(pS_1) = Rs(pS_1) = p$ as no move is available. Thus $m_1 = 1$.

\noindent If $n=3$, any move of a player removes all the three vertices of $S_3$. $Ls$ and $Rs$ can thus be computed as follows: 
\begin{align*}
&Ls(2pS_3) = 3 + Rs\big((2p-1)S_3\big) = 3 - 3 + Ls\big((2p-2)S_3\big) = \dots = 0 \\
&Ls\big((2p+1)S_3\big) = 3 + Rs\big(2pS_3\big) = 3 - 3 + Ls\big((2p-1)S_3\big) = \dots = 3 \\
\end{align*}
\noindent This proves that $0 \le Ls(pS_3) \le 3$ for all $p$. Consequently, $m_3 = \underset{p \to \infty}{\lim} \frac{Ls(pS_3)}{p} = 0$.

\noindent If $n=5$, by  Example~\ref{exemple P5}, we have $4S_5 = 4$. By denoting $p = 4q + r$ with $0 \le r <4$, we have $Ls(pS_5) = Ls(rS_5) + 4q$. Thus, since $0 \le r < 4$, $Ls(rS_5)$ is bounded, $\underset{p \to \infty}{\lim} \frac{Ls(pS_{5})}{p} = 1$.

\vspace{.2cm}
Suppose $n\ge 7$. We first prove that $Ls(pS_n) \le p + 8$. Assuming that Left starts, we denote by $p_k$ the number of segments of odd length in the resulting graph after the $k$-th move, and by $s_k$ the relative score, i.e. the difference between the number of vertices taken by Left and the number of vertices taken by Right after the $k$-th move. We will define a strategy for Right that has the following invariant property: for any $k \ge 0$, if there exists at least one segment of length four or greater, we have $p_{2k} + s_{2k} \le p$. 
Since we have $p_0 = p$ and $s_0 = 0$, the invariant is true before the first move. Note that we will only consider $p_{2k}$ and $s_{2k}$ to focus on the state of the game after each move of Right.

\noindent The strategy of Right is the following: while there exists one segment of size at least four, Right considers a segment of highest length (if there are several, consider any of them). If it has odd length, Right can, by playing a $3$-move on it on an extremity, transform it into a segment of even length; if it has even length, Right can, by playing a $4$-move on it (which is always available on an even length segment of size greater than four) take four vertices and leave a segment of the same parity. Therefore, after the $k$-th move of Right, we have $s_{2k} + p_{2k} = s_{2k-1} + p_{2k-1} -4$. Note that with such a strategy, Right does not create new odd segment.


\noindent Now, we consider the variation of the sequences $s_k$ and $p_k$ after the $(k+1)$-th move of Left.

\begin{itemize}
    \item If the move of Left is a $5$-move, necessarily, it is played on a $S_5$ that is removed. We have $p_{2k+1} = p_{2k}-1$ and $s_{2k+1} = s_{2k} + 5$, thus $p_{2k+1} + s_{2k+1} = p_{2k} + s_{2k} +4$.
    
    \item If the move of Left is a $4$-move, it does not change the number of segments of odd length. We have $p_{2k+1} = p_{2k}$ and $s_{2k+1} = s_{2k} + 4$, and thus $p_{2k+1} + s_{2k+1} = p_{2k} + s_{2k} +4$.

    \item If the move of Left is a $3$-move, and if it is in an even segment, it creates one odd segment and one even segment. If it is in an odd segment, it creates two odd segments, unless it is on a $S_3$. In both cases, the number of odd segments is increased by at most one. We have $p_{2k+1} \le p_{2k} +1$ and $s_{2k+1} = s_{2k} + 3$, thus $p_{2k+1} + s_{2k+1} \le p_{2k} + s_{2k} +4$. Note that since the odd segments are positive (i.e. $n\geq 0$) at the beginning of the game, and since Right's strategy does not create new odd segment, Left cannot create negative odd segment when playing. 
    
    \item If the last move of Left is a $2$-move, it does not change the parity of the number of odd segments. We have $p_{2k+1} = p_{2k}$ and $s_{2k+1} = s_{2k} + 2$, thus $p_{2k+1} + s_{2k+1} = p_{2k} + s_{2k} +2$.
\end{itemize}

\noindent Finally, we always have $p_{2k+1} + s_{2k+1} \le p_{2k} + s_{2k} +4$. Thus, as $p_{2k} + s_{2k} = p_{2k-1} + s_{2k-1} -4$, the invariant property is satisfied until all segments have length three or less. In the final part of the game when all segments have length three or less, as any move in a segment removes it fully, a greedy strategy ensures that Left takes at most three vertices more than Right. As the largest available value for $s_k$ is $p + 5$ (i.e., after a $5$-move of Left, according to the invariant), we obtain $Ls(pS_n) \le p + 5+3 = p+8$.

\vspace{.2cm}

It remains to prove that $Ls(pS_n) \ge 0$. We consider a strategy for Left that preserves the following invariant: at any moment of the game, Left has taken at least the same number of vertices than Right. Moreover, if it is Right's turn, there is no segment $S_k$ with $k<0$, i.e. all the remaining segments have at least one extremity black, and Left has taken at least four vertices more then Right.

\noindent We define here the first part of the strategy for Left. Until all segments have size 3 or less, Left considers the following strategy:

\begin{itemize}
    \item First, Left takes four vertices on a segment. This move is available as we assumed $n \ge 7$.
    
    \item If the last move of Right does not create a segment with both extremities white, if at least one segment has size four or more, Left plays a $4$-move on it (it might be a $5$-move if the segment has length 5). In this way, Left has taken at least four vertices and Right at most four, as $5$-moves are only available for Right on $S_{-5}$ and the invariant ensures that there is no segment of this form when it is Right's turn. Therefore, the invariant is satisfied as Left does not create a $S_k$ with $k<0$.
    
    \item If the last move of Right creates a segment with both extremities white, by induction, such a segment did not exist before he made this move. Thus, it has been created by a $3$-move on a segment with extremities of different colors. Therefore, this move has created a segment with two white extremities and an even segment (eventually of size zero). Left answers by playing a $3$-move on an extreme black vertex of the segment which has two white extremities. Playing so, Left takes also three vertices and this segment now has a black extremity. Once again, the invariant is still satisfied.
\end{itemize}

 \noindent This strategy guarantees that the invariant is satisfied until all segments have size $3$ or less. At this moment of the game, if it is Left's turn, by a greedy strategy, Left can ensure to take at least the same number of vertices than Right for the rest of the game. Thus, as she already has at least the same number of vertices that Right, Left is not losing.
 
 \noindent If it is Right's turn, by the invariant property, Left has taken at least four vertices more than Right (the four vertices played at the beginning). Thus, regardless the move that Right does, Left will get more vertices than Right as any segment has size three or less. This strategy ensures $Ls(pS_n) \ge 0$.
 
 By gluing the two above results, we get $0 \le Ls(pS_n) \le p + 8$, implying $ 0 \le \frac{Ls(pS_n)}{p} \le 1 + \frac{8}{p}$. By taking the limit, as we know from Theorem \ref{thmmoy} that $m_n$ exists, we have $0 \le m_n \le 1$. 
 
\end{proof}

The above result about the mean is very useful to estimate the score of a sum of segments according to Corollary~\ref{coro:sumtemp}. Yet, it remains to have an upper bound on the temperature of a segment. The following theorem provides such a bound.

\begin{theorem}
Let $n \in \Z$. We have $\si(S_n) \le 4$.
\label{thmtempseg}
\end{theorem}

\begin{proof}
Up to switching Right and Left, we can assume without loss of generality that $n \ge 0$. We prove the result by induction on $n$. If $n \le 4$, the result is straightforward as the temperature cannot be larger than the maximum score on the game. If $n = 5$, one can prove that $\si(S_n) = 4$ according to Figure~\ref{fig:treecooled}.

Assume now that $n \ge 6$.  By Definition \ref{def:temp} of the temperature, to prove that $\si(S_n)\le 4$, it is enough to prove that for $t\geq 4$, $Ls((\widetilde{S_n})_t)\leq Rs((\widetilde{S_n})_t)$.
Let $t\ge 4$. We recall that  $(\widetilde{S_n})_t= \langle (S_n^L)_t-t|(S_n^R)_t+t \rangle$. Consider a Left option $g^L$ and a Right option $g^R$ of $S_n$. 


A move in $S_n$ splits the segment into at most two parts, thus we can write $g^L = \al +  S_{i^L_1} + S_{i^L_2}$, and $g^R = - \be + S_{i^R_1} + S_{i^R_2}$ with $\al, \be >0$ the scores obtained by each move.

By induction, we have $\si(S_{i^L_1}), \si(S_{i^L_2}), \si(S_{i^R_1}),  \si(S_{i^R_2}) \le 4$. Thus, $g^L_t$ and $g^R_t$ are frozen for $t\geq 4$ and, by Proposition \ref{prop:sim}, equal to their means. In particular, since the mean is linear, we have $Rs(g^L_t)=\al+m(S_{i^L_1})+m(S_{i^L_2})$ and
$Ls(g^R_t)=-\be+m(S_{i^R_1})+m(S_{i^R_2})$.

Assume first that $n$ is odd. Right can only play $3$-moves from $S_n$, thus $\beta=3$ and both $S_{i^R_1}$ and $S_i{^R_2}$ are of even length. Therefore, by Theorem~\ref{thm moy segment}, $m(S_{i^R_1}) = m(S_{i^R_2}) = 0$ and we have: 

$$Ls(g^R_t) + t = -3 + t + m(S_{i^R_1}) + m(S_{i^R_2}) = -3 + t.$$

If Left made a $4$-move, up to exchange $S_{i^L_1}$ and $S_{i^L_2}$, we can suppose $S_{i^L_2} = 0$. By Theorem~\ref{thm moy segment}, $m(S_{i^L_1}) \le 1$ and thus: 

$$Rs(g^L_t) - t = 4 - t + m(S_{i^L_1}) \le 4 - t +1 = 5 - t.$$

If Left made a $3$-move,  $m(S_{i^L_1}), m(S_{i^L_2}) \le 1$ by Theorem~\ref{thm moy segment}, and we have: 

$$Rs(g^L_t) - t = 3 - t + m(S_{i^L_1}) + m(S_{i^L_2}) \le 3 - t + 1 + 1 = 5 - t.$$

Finally, we obtain, for any options $g^L$ and $g^R$:
\begin{align*}
   Rs(g^L_t) - Ls(g^R_t) & \le 5 - t - ( -3 + t) \\
   & \le 8 - 2t \\
   & \le 0
\end{align*}

In particular, this is true if $g^L_t-t$ and $g^R_t+t$ are the best options of $\widetilde{(S_n)}_t$. Then we have: 

    $$Ls(\widetilde{(S_n)}_t) - Rs(\widetilde{(S_n)}_t) = Rs(g^L_t) - Ls(g^R_t)  \le 0.
$$

Assume now that $n$ is even. If Left made a $3$-move, one of the remaining path is even and the other is odd. Without loss of generality, we can assume that $S_{i^L_2}$ has even length. By Theorem \ref{thm moy segment}, $m(S^L_{i_2}) = 0$ and $|m(S^L_{i_1})| \le 1$. This implies that
    $$Rs(g^L_t) - t = 3 + m(S^L_{i_1}) + m(S^L_{i_2}) - t \le 4 - t.$$  
  
Otherwise, Left made a $4$-move and only one path is remaining, and has even length . Thus we also have:
   $$Rs(g^L_t) - t = 4 + m(S_{n-4})-t\le 4 - t.$$
  
    Symmetrically, $Ls(G^R_t) + t \leq -4+t$. 
     Therefore, if $t \ge 4$, we have for any options $g^L$ and $g^R$, 
    $Rs(g^L_t) - Ls(g^R_t) \le 8 - 2t \leq 0$ and so, as before, $Ls(\widetilde{(S_n)}_t) - Rs(\widetilde{(S_n)}_t) \leq 0$.

Finally, regardless if $n$ is odd or even, we have  $Ls(\widetilde{(S_n)}_t) - Rs(\widetilde{(S_n)}_t)\leq 0$ for any $t\geq 4$. In particular, this means that $\si(S_n) \le 4$.
\end{proof}

We now have all the material to give an estimate of the score on a sum of segments.

\begin{corollaire}\label{cor:sum}
Let $n_1, \dots, n_p \in \Z$, with exactly $k$ of the $n_i$ of odd value. We have $$-k - 4 \le Rs(S_{n_1} + \dots + S_{n_p}) \le Ls(S_{n_1} + \dots + S_{n_p}) \le k + 4$$.
\end{corollaire}

\begin{proof}
This result is directly deduced from Theorem~\ref{thmtempseg}, Theorem~\ref{thm moy segment} and Corollary \ref{coro:sumtemp}.
\end{proof}

Note that the bounds of Corollary~\ref{cor:sum} are tight. Consider for example a sum of $4n+1$ segments $S_5$. Since $4 S_5 = 4$ according to Example~\ref{exemple P5}, we have $Ls( (4n+1) S_5) = 4n + Ls(S_5) = 4n+5$. 

\subsection{Additional results about the periodicity conjecture of segments}

In order to explore the validity of Conjecture~\ref{conj:periodicity}, more values are required. Recall that in \cite{duchene}, only the first $80$ values of the score has been computed. The computation of the scores of segments can be boosted by using equivalences between sums. Recall that in a sum of games, any component can be replaced by another component equivalent to it.

The following result shows that two consecutive segments of even size (starting by the one that is a multiple of $4$) can be replaced by the segment $S_2$ in a sum. 

\begin{theorem}
\label{thm4k-4k2}
Let $k \in \N$, we have $S_{4k} + S_{4k+2} = S_2$. 
\end{theorem}

\begin{proof}
According to Lemma~\ref{lemma equiv}, it suffices to show that $Ls(S_{4k} + S_{4k+2} + S_2) = 0$.

We will prove by induction on $k$ that for any $k\in \N$, we have :

\begin{align*}
Ls(S_{-(4k+5)} + S_{4k+3} + S_2) = 0 \\
Ls( S_{-(4k+3)} + S_{4k+1} + S_2) = 0 \\
Ls(S_{4k} + S_{4k+2} + S_2) = 0
\end{align*}

First, we give strategies for Left to ensure the following inequalities:

\begin{align*}
Ls(S_{-(4k+5)} + S_{4k+3} + S_2) \ge 0 \\
Ls( S_{4k+3} + S_{-(4k+1)} + S_2) \ge 0 \\
Ls(S_{4k} + S_{4k+2} + S_2) \ge 0
\end{align*}

\noindent In $S_{4k} + S_{4k+2} + S_2$, Left can play a move on a extremal vertex in $S_{4k+2}$ that takes exactly two vertices. After that, it is Right's turn and the game is $S_{4k} + S_{4k} + S_2$ that is equivalent to $S_2$. Therefore, Right also scores two points and we have $Ls(S_{4k+2} + S_{4k} + S_2)\ge 0$.

\noindent In $S_{-(4k+5)} + S_{4k+3} + S_2$ , Left can first play $S_2$ and scores two points. Then, by Theorem~\ref{thm:givingvertices} (third inequality), imagine that Left gives the second vertex of $S_{-(4k+5)}$ to Right. After this move, the game is equivalent to $S_{-(4k+3)} + S_{4k+3} + S_{-1} = -1$. Once again, we have $Ls( S_{-(4k+5)} + S_{4k+3} + S_2) \ge 0$. 

\noindent In $S_{4k+3} + S_{-(4k+1)} + S_2$ , Left can first play an extremal move on $S_{4k+3}$ that takes exactly two vertices. After that, it is Right's turn and the game is $S_{4k+1} + S_{-(4k+1)} + S_2 = S_2$. Therefore, Right scores two points on it and we have $Ls(S_{4k+3} + S_{-(4k+1)} + S_2)\ge 0$.

\vspace{.2cm}

Now, we give strategies for Right that ensure the other inequalities (i.e. Right is the second player and is not losing).

\noindent If $k = 0$ or $k = 1$, the induction hypothesis is true by checking by hand all the possible options. Now assume $k \ge 2$ and that the induction hypothesis is true for all $0\leq k'<k$.

\noindent First, remark that no 5-move is available in any of these games for Left as there is no $S_5$ component.

\noindent Secondly, if Left plays a $4$-move in any of these games, Right can answer by a $4$-move in the other component that is not $S_2$ and then obtains an instance that is strictly smaller that the initial one and where he will not be losing by induction.

\noindent Thirdly, since $2$-moves are dominated, one can suppose, without loss of generality, that Left never plays one of them.

\noindent It thus remains to consider which strategy can be used by Right when Left starts with a $3$-move.

\noindent If the game is $S_{4k} + S_{4k+2} + S_2$:
\begin{itemize}

    \item If Left plays a $3$-move in $S_{4k}$ :
        \begin{itemize}
            \item if this move creates two instances $S_{4m+1} + S_{4n}$ (with $0<n,m<k$), Right can answer by playing a $3$-move in $S_{4k+2}$ that creates $S_{-(4m+1)} + S_{4n + 2}$. The resulting game is $S_{4m+1} + S_{4n} + S_{-(4m+1)} + S_{4n + 2} + S_2 = S_{4n} + S_{4n+2} + S_2$ that is not losing for Right by induction hypothesis.
            
            \item if this move creates two instances $S_{4m+3} + S_{4n+2}$, Right can answer by playing a $3$-move in $S_{4k+2}$ that creates $S_{-(4m+5)} + S_{4n + 2}$. The resulting game becomes $S_{-(4m+5)} + S_{4n+2} + S_{4m+3} + S_{4n + 2} + S_2 = S_{4m+3} + S_{-(4m+5)} + S_2$ that is also not losing for Right by induction hypothesis.
        \end{itemize}

    \item If Left plays a $3$-move in $S_{4k+2}$:
    \begin{itemize}
        \item if she creates two instances $S_{4m+2} + S_{4n+1}$, Right can answer by playing a $3$-move in $S_{4k}$ that creates $S_{4m} + S_{-(4n+1)}$.
        The game becomes $S_{4m+2} + S_{4n+1} + S_{4m} + S_{-(4n+1)} + S_2 = S_{4m} + S_{4m+2} + S_2$ and we conclude as previously by induction.
        
        \item if she creates two instances $S_{4m} + S_{4n+3}$, Right can answer by playing a $3$-move in $S_{4k}$ that creates $S_{4m} + S_{-(4n+1)}$. The game becomes $S_{4m} + S_{4n+3} + S_{4m} + S_{-(4n+1)} + S_2 = S_{4n+3} + S_{-(4n+1)} + S_2$ and we can conclude by induction.
        
    \end{itemize}
\end{itemize}

\noindent If the game is $S_{4k+3} + S_{-(4k+5)} + S_2$:

\begin{itemize}
    \item if Left makes a $3$-move in $S_{4k+3}$, it necessarily creates two instances $S_{4m+1}$ and $S_{4n+3}$. Right can answer by playing a $3$-move in $S_{-(4k+5)}$ that creates $S_{-(4m+1)}$ and $S_{-(4n+5)}$. The resulting game is $S_{4m+1} + S_{4n+3} + S_{-(4m+1)} + S_{-(4n+5)} + S_2 = S_{4n+3} + S_{-(4n+5)} + S_2$, and thus not losing for Right by induction hypothesis.
    \item if Left makes a $3$-move in $S_{-(4k+5)}$, it creates two instances $S_{4m}$ and $S_{4n+2}$. Right can answer by playing a $3$-move in $S_{-(4k+3)}$ that creates $S_{4m}$ and $S_{4n}$. The game becomes $S_{4m} + S_{4n +2} + S_{4m} + S_{4n} + S_2 = S_{4n} + S_{4n + 2} + S_2$, and we can conclude by induction.
\end{itemize}

\noindent If the game is $S_{4k+3} + S_{-(4k+1)}$:

\begin{itemize}
    \item if Left makes a $3$-move in $S_{4k+3}$, she creates two instances $S_{4m+1}$ and $S_{4n+3}$. Right can answer by playing a $3$-move in $S_{-(4k+1)}$ that creates $S_{-(4m+1)}$ and $S_{-(4n+1)}$. The resulting game is $S_{4m+1} + S_{4n+3} + S_{-(4m+1)} + S_{-(4n+1)} + S_2 = S_{4n+3} + S_{-(4n+1)} + S_2$, and hence not losing for Right by induction.

    \item if Left makes a $3$-move in $S_{-(4k+1)}$, she creates two instances $S_{4m+2}$ and $S_{4n}$. Right can answer by playing a $3$-move in $S_{4k+3}$ that creates $S_{4m+2}$ and $S_{4n+2}$. The resulting game becomes $S_{4m+2} + S_{4n} + S_{4m+2} + S_{4n+2} + S_2 = S_{4n} + S_{4n + 2} + S_2$ and we can conclude by induction.
\end{itemize}

\end{proof}

The above result yields a very significant way to make the computation of the values faster. In fact, as our program is based on a dynamic programming approach, all sums of segments need to be stored. Thanks to this result, all sums containing a $S_{4k+2}$ can be simplified by replacing it with a $S_{4k} + S_2$.\\

Table~\ref{table120} expands the result presented in~\cite{duchene} by giving the first $120$ values of the score (instead of $80$). A careful analysis of this table suggests that Conjecture~\ref{conj:periodicity} holds, with a period of $40$ and a preperiod of $30$. In addition, one can remark that after the preperiod, the structure of the values is very regular (almost of period $8$), except for $n=37,77,117$ that correspond to the rare values satisfying $Rs(S_n)=-5$. \\

This new table is however not sufficient to prove the ultimate periodicity of the structure, as the conjectured period of $40$ does not appear frequently enough to be used in a general proof. In addition, the segments $S_n$ and $S_{n+40}$ are not equivalent for almost all values of (up to $n = 40$), which is another pitfall to prove the conjecture. A first step towards it would be to prove that the values of $Ls(S_n)$ are equal to $2$ for all even $n\geq 32$.

\begin{figure}
    \centering
    
    \scalebox{.75}{ 
      \begin{tabular}{|c||r|r|r|r|r|r|r|r|r|r|r|r|r|r|r|r|r|r|r|r|r|}
     \hline
n & 001 & 002 & 003 & 004 & 005 & 006 & 007 & 008 & 009 & 010 & 011 & 012 & 013 & 014 & 015 & 016 & 017 & 018 & 019 & 020\\ 
 \hline
Ls($S_n$) & 1 & 2 & 3 & 4 & 5 & 2 & 1 & 2 & 3 & 2 & 1 & 2 & 3 & 4 & 3 & 2 & 3 & 2 & 3 & 2\\ 
Rs($S_n$) & 1 & -2 & -3 & -4 & -1 & -2 & -3 & -2 & -1 & -2 & -3 & -2 & -1 & -4 & -3 & -2 & -1 & -2 & -3 & -2\\
\hline
   \end{tabular}   }
   
   \scalebox{.75}{
   
      \begin{tabular}{|c||r|r|r|r|r|r|r|r|r|r|r|r|r|r|r|r|r|r|r|r|r|}
     \hline
n & 021 & 022 & 023 & 024 & 025 & 026 & 027 & 028 & 029 & 030 & 031 & 032 & 033 & 034 & 035 & 036 & 037 & 038 & 039 & 040\\ 
 \hline
Ls($S_n$) & 5 & 4 & 3 & 2 & 3 & 2 & 3 & 2 & 5 & 4 & 3 & 2 & 3 & 2 & 3 & 2 & 5 & 2 & 3 & 2\\ 
Rs($S_n$) & -1 & -4 & -3 & -2 & -1 & -2 & -3 & -2 & -1 & -4 & -3 & -2 & -1 & -2 & -3 & -2 & -1 & -2 & -3 & -2\\
\hline
   \end{tabular}   
   }
   
   \scalebox{.75}{

           \begin{tabular}{|c||r|r|r|r|r|r|r|r|r|r|r|r|r|r|r|r|r|r|r|r|r|}
     \hline
n &  041 &  042 &  043 &  044 &  045 &  046 &  047 &  048 &  049 &  050 &  051 &  052 &  053 &  054 &  055 &  056 &  057 &  058 &  059 &  060\\ 
 \hline
Ls($S_n$) &  3 & 2 & 3 & 2 & 3 & 2 & 3 & 2 & 3 &  2 & 3 & 2 & 3 & 2 &  3 & 2 & 3 & 2 & 3 & 2\\ 
Rs($S_n$) &  -1 &  -2 &  -3 &  -2 &  -1 &  -2 &  -3 &  -2 &  -1 &  -2 &  -3 &  -2 &  -1 &  -2 &  -3 &  -2 &  -1 &  -2 &  -3 &  -2\\
\hline

   \end{tabular}  }
   
   \scalebox{.75}{

              \begin{tabular}{|c||r|r|r|r|r|r|r|r|r|r|r|r|r|r|r|r|r|r|r|r|r|}
     \hline
n &  061 &  062 &  063 &  064 &  065 &  066 &  067 &  068 &  069 &  070 &  071 &  072 &  073 &  074 &  075 &  076 & 077 & 078 & 079 & 080\\ 
 \hline
 
Ls($S_n$) & 3 & 2 & 3 & 2 & 3 & 2 & 3 & 2 & 3 & 2 & 3 & 2 & 3 & 2 & 3 & 2 & 5 & 2 & 3 & 2\\ 
Rs($S_n$) & -1 & -2 & -3 & -2 & -1 & -2 & -3 & -2 & -1 & -2 & -3 & -2 & -1 & -2 & -3 & -2 & -1 & -2 & -3 & -2\\
\hline
   \end{tabular} 
   
   }
   
   \scalebox{.75}{
   
               \begin{tabular}{|c||r|r|r|r|r|r|r|r|r|r|r|r|r|r|r|r|r|r|r|r|r|}
     \hline
n & 081 & 082 & 083 & 084 & 085 & 086 & 087 & 088 & 089 & 090 & 091 & 092 & 093 & 094 & 095 & 096 & 097 & 098 & 099 & 100\\ 
 \hline
Ls($S_n$) & 3 & 2 & 3 & 2 & 3 & 2 & 3 & 2 & 3 & 2 & 3 & 2 & 3 & 2 & 3 & 2 & 3 & 2 & 3 & 2\\ 
Rs($S_n$) & -1 & -2 & -3 & -2 & -1 & -2 & -3 & -2 & -1 & -2 & -3 & -2 & -1 & -2 & -3 & -2 & -1 & -2 & -3 & -2\\
\hline
   \end{tabular} 
   }
   
   \scalebox{.75}{
   
              \begin{tabular}{|c||r|r|r|r|r|r|r|r|r|r|r|r|r|r|r|r|r|r|r|r|r|}
     \hline
n & 101 & 102 & 103 & 104 & 105 & 106 & 107 & 108 & 109 & 110 & 111 & 112 & 113 & 114 & 115 & 116 & 117 & 118 & 119 & 120\\
\hline
Ls($S_n$) & 3 & 2 & 3 & 2 & 3 & 2 & 3 & 2 & 3 & 2 & 3 & 2 & 3 & 2 & 3 & 2 & 5 & 2 & 3 & 2\\
Rs($S_n$) & -1 & -2 & -3 & -2 & -1 & -2 & -3 & -2 & -1 & -2 & -3 & -2 & -1 & -2 & -3 & -2 & -1 & -2 & -3 & -2\\
     \hline
   \end{tabular}

   }
    
    \caption{New values on Segments}
    \label{table120}
\end{figure}

\section{Grids}\label{sec:grids}

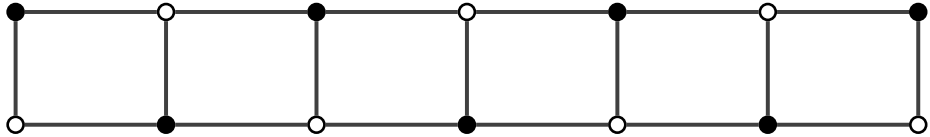
\begin{figure}[ht]
    \centering
\begin{tikzpicture}

\draw (2,1.5) node[vertexB](3){} ;

    \draw (2,0) node[vertexW](4){} ;

    \draw (4,0) node[vertexB](5){} ;
  \draw (4,1.5) node[vertexW](6){} ;

  \draw (6,1.5) node[vertexB](7){} ;

  \draw (6,0) node[vertexW](8){} ;

  \draw (8,0) node[vertexB](9){} ;

  \draw (8,1.5) node[vertexW](10){} ;

  \draw (10,1.5) node[vertexB](11){} ;

  \draw (10,0) node[vertexW](12){} ;

  \draw (12,0) node[vertexB](13){} ;

  \draw (12,1.5) node[vertexW](14){} ;

  \draw (14,1.5) node[vertexB](15){} ;

  \draw (14,0) node[vertexW](16){} ;

  \Edge[](4)(3) 
  \Edge[](4)(5)
  \Edge[](6)(3)
  \Edge[](6)(5)
  \Edge[](6)(7)
  \Edge[](8)(5)
  \Edge[](8)(7)
  \Edge[](8)(9)
  \Edge[](10)(7)
  \Edge[](10)(9)
  \Edge[](10)(11)
  \Edge[](12)(9)
  \Edge[](12)(11)
  \Edge[](12)(13)
  \Edge[](14)(11)
  \Edge[](14)(13)
  \Edge[](14)(15)
  \Edge[](16)(13)
  \Edge[](16)(15)

\end{tikzpicture}

    \caption{Grid $G_{2,7}$}
    \label{fig:G27}
\end{figure}

As a natural extension to paths, we consider in this section rectangular grids.
We denote by $G_{n,m}$ the grid with $n$ rows and $m$ columns, with alternated black and white vertices and a black vertex in the top left corner. See Figure \ref{fig:G27} for an illustration of $G_{2,7}$.

Despite their symmetry, there is no BW-automorphism in rectangular grids. Indeed, the vertices close to the center of the grid are necessarily sent to each other by any automorphism, and thus could not be sufficiently distant. The results of Section \ref{sec:symetrie} cannot therefore be applied.
Moreover, it seems that there are neither draws nor games favourable to the second player in rectangular grids. Actually, all the computations made so far lead to a victory of the first player. This induces the following conjecture.

\begin{conjecture}\label{conj:grids}
Let $n,m \geq 2$. We have $Ls(G_{n,m}) > 0 > Rs(G_{n,m})$. In particular, the first player wins and there is no draw.
\end{conjecture}

Note that if $n$ or $m$ is even, the grid is symmetric by exchanging the roles of Left and Right. Since the game is nonzugzwang, $Ls(G_{n,m}) = - Rs(G_{n,m})$. Moreover, as $Ls(G_{n,m}) \ge Rs(G_{n,m})$, it implies directly that $Ls(G_{n,m}) \ge 0 \ge Rs(G_{n,m})$. However, if $n$ and $m$ are odd, this argument cannot be used.

In the rest of the section, we give some partial results supporting the conjecture for rectangular grids with two and three rows.
 
\subsection{Grids with two rows}
In this subsection, we prove Conjecture \ref{conj:grids} on grids with only two rows. As said before, since the number of rows is even, we just need to consider the Left score. We actually prove a stronger result when the number of columns is odd, by giving the exact value of $Ls(G_{2,m})$.

\begin{theorem}\label{thm:2rowsodd}
Let $m\geq 3$ be an odd integer. We have $Ls(G_{2,m}) = -Rs(G_{2,m}) = 4$ if $m\geq 5$ and $Ls(G_{2,3})=-Rs(G_{2,3})=6$.
\end{theorem}

For even $m$, we just manage to prove a lower bound on the score.

\begin{lemma}\label{lem:2rowseven}
  Let $m\geq 4$ be an even integer. We have $Ls(G_{2,m})\geq 2$.
\end{lemma}

Before proving Theorem \ref{thm:2rowsodd} and Lemma \ref{lem:2rowseven}, we introduce some notations that will be used in both proofs. As each column contains at most one  white and one black vertex, we will say playing $i$ for playing the vertex of the current player's color in column $i$.
Consider a connected component $G$ at some point of a game played in a grid with two rows. Note that $G$ can only have columns with one vertex at its extremities. Thus the full columns are consecutive. We encode $G$ as a triplet $(x,n,y)$ with  $x,y \in \{0,-1,1\}$ and $n$ an integer. The integer $n$ corresponds to the number of consecutive full columns, $x$ represents the first column of $G$ and $y$ the last one with the following convention: $x=0$ ($y=0$ resp.) if the column is full, $x=1$ if there is a black vertex and $x=-1$ if there is a white vertex. In Figure \ref{fig:Ls(G)<=4}, the two grids have respectively forms $(-1,2-1)$ and $(1,4,1)$.  Note that the form of a partial grid completely characterizes it and that $(x,n,y)$ and $(y,n,x)$ correspond to the same graph.


\begin{proof}[Proof of Theorem \ref{thm:2rowsodd}] 
For $m=3$, the first player can take the whole grid in one move, which implies that $Ls(G_{2,3})=-Rs(G_{2,3})=6$.

Let $m \geq 5$ be an odd integer.
We first prove that $Ls(G) \ge 4$. Indeed, Left can play $\frac{m+1}{2}$, i.e in the middle column.
This way, she takes four vertices and cuts the game in two parts of form $(0,\frac{n-3}{2},1)$ (see Figure~\ref{fig:odd}).
Then, each time Right decides to play on column $i$, Left plays symmetrically on column $m-i+1$. If Right plays
$\frac{m+1}{2}-2$ (respectively  $\frac{m+1}{2}+2$) and removes the single black vertex of column $\frac{m+1}{2}-1$ (resp. $\frac{m+1}{2}+1$), the symmetric move of Left will remove the symmetric vertices of Right's move plus the other single black vertex. In all the other cases, Left's move removes exactly the symmetric vertices (considering a central symmetry) that Right has taken.
 This way, Left can finish the game without loosing any vertex and $Ls(G_n) \ge 4$.

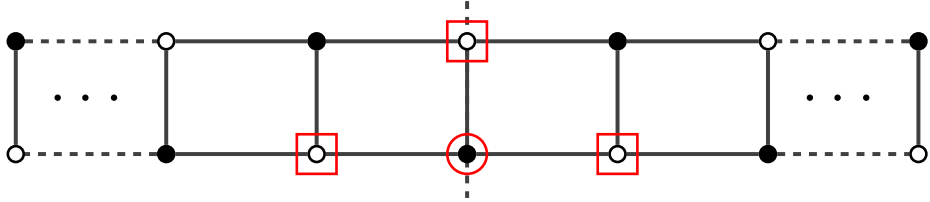
\begin{figure}[ht]
    \centering

\begin{tikzpicture}

\Vertex[x=1,y=-.75, style = {color = white}, label = \dots, size = .5, fontscale = 3]{16}

\Vertex[x=11, y=-.75, style = {color = white}, label = \dots, size = .5, fontscale = 3]{16}

  \draw (0,0) node[vertexB](1){} ;

  \draw (0,-1.5) node[vertexW](2){} ;

  \draw (2,-1.5) node[vertexB](3){} ;

  \draw (2,0) node[vertexW](4){} ;

  \draw (4,0) node[vertexB](5){} ;

  \draw (4,-1.5) node[vertexW](6){} ;
  \draw (4,-1.5) node[redsquare](n){} ;

  \draw (6,-1.5) node[vertexB](7){} ;
  \draw (6,-1.5) node[redcircle](n){} ;

  \draw (6,0) node[vertexW](8){} ;
  \draw (6,0) node[redsquare](n){} ;

  \draw (8,0) node[vertexB](9){} ;

  \draw (8,-1.5) node[vertexW](10){} ;
  \draw (8,-1.5) node[redsquare](n){} ;

  \draw (10,-1.5) node[vertexB](11){} ;

  \draw (10,0) node[vertexW](12){} ;

  \draw (12,0) node[vertexB](13){} ;

  \draw (12,-1.5) node[vertexW](14){} ;

  \Vertex[x = 6, y = -2.2, size = .01, style = {color = white}]{20}

  \Vertex[x = 6, y = .7, size = .01, style = {color = white}]{21}

  \Edge[, style = dashed](1)(4)
  \Edge[](1)(2)
  \Edge[, style = dashed](3)(2)
  \Edge[](3)(4)
  \Edge[](5)(4)
  \Edge[](3)(6)
  \Edge[](5)(6)
  \Edge[](7)(6)
  \Edge[](5)(8)
  \Edge[](7)(8)
  \Edge[](9)(8)
  \Edge[](7)(10)
  \Edge[](9)(10)
  \Edge[](11)(10)
  \Edge[](9)(12)
  \Edge[](11)(12)
  \Edge[style = dashed](13)(12)
  \Edge[style = dashed](11)(14)
  \Edge[](13)(14)
  \Edge[style = dashed](20)(21)

\end{tikzpicture}

    \caption{An optimal strategy for Left when the number of columns is odd, is to take the central black vertex (marked with a red circle) and then playing the symmetric column. This strategy gives a final score of 4.}
    \label{fig:odd}
\end{figure}

We now prove that $Ls(G_{2,m}) \le 4$ by induction. First, we  prove that if $G$ has the form $G_1 + G_2$ with $G_1 = (x,n,y)$ and $G_2 = (-x,n+2,-y)$, then $Ls(G) \le 4$ (see Figure \ref{fig:Ls(G)<=4}).

\begin{figure}
    \centering
    
\begin{tikzpicture}

  \draw (4,0) node[vertexB](-2){} ;

  \draw (2,0) node[vertexW](-1){} ;

  \draw (2,1.5) node[vertexB](0){} ;

  \draw (0,0) node[vertexB](1){} ;

  \draw (0,1.5) node[vertexW](2){} ;

  \draw (-2,1.5) node[vertexB](3){} ;

  \draw (-2,0) node[vertexW](4){} ;

  \draw (-4,0) node[vertexB](5){} ;
  
  \draw (-4,1.5) node[vertexW](6){} ;
  
  \draw (-6,1.5) node[vertexB](7){} ;

  \draw (-8,1.5) node[vertexW](10){} ;

  \draw (-10,1.5) node[vertexB](11){} ;

  \draw (-10,0) node[vertexW](12){} ;

  \draw (-12,0) node[vertexB](13){} ;

  \draw (-12,1.5) node[vertexW](14){} ;

  \draw (-14,0) node[vertexW](15){} ;

  \Edge[](-1)(1)
  \Edge[](-1)(-2)
  \Edge[](0)(-1)
  \Edge[](0)(2)
  \Edge[](1)(4)
  \Edge[](1)(2)
  \Edge[](3)(2)
  \Edge[](3)(4)
  \Edge[](5)(4)
  \Edge[](3)(6)
  \Edge[](5)(6)
  \Edge[](7)(6)

  \Edge[](11)(10)
  \Edge[](11)(12)
  \Edge[](13)(12)
  \Edge[](11)(14)
  \Edge[](13)(14)
  \Edge[](13)(15)

\end{tikzpicture}

    \caption{In $(-1,2,-1)+(1,4,1)$, Left can score at most $4$.}
    \label{fig:Ls(G)<=4}
\end{figure}
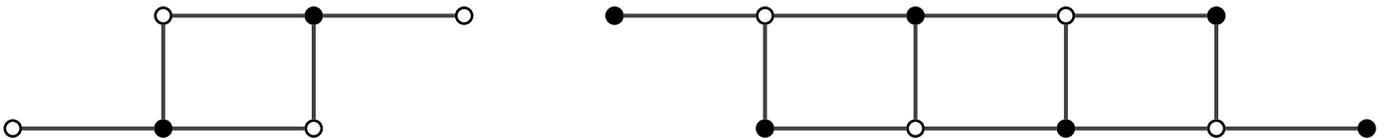

We prove this result by induction on $n$ the number of full columns of $G_1$ ($G_2$ has $n+2$ full columns). 
Note that $G_2$ has four more vertices than $G_1$. Consider the first move of Left. Assume first that Left did not take a full component and that Left plays in $G_1$. By symmetry, we can assume Left has play in some column $i$ with $i\leq \frac{n+1}{2}$ (starting to count from the first full column).
The whole graph after this move can be written as
$R+(1,n-i-1,y) +G_2$ with $R$ that might be empty (and necessarily, $i+1\leq n$).
Then Right can answer by playing in column $i$ of $G_2$. He takes the same number of vertices than Left. The remaining graph after these two moves has the form $R+(1,n-i-1,y)+R'+(-1,n-i+1,-y)$ where $R$ and $R'$ are opposite (and might be empty). Hence, $R$ and $R'$ cancelled each other. By induction, Left can score at most $4$ in $(1,n-i-1,y)+(-1,n-i+1,-y)$, and thus can score at most $4$ in total.

Assume now that Left plays in $G_2$ but did not take the whole graph. Again, by symmetry, we can assume Left has play in some column $i$ with $i\leq \frac{n+3}{2}$ and that it remains $G_1+R'+(1,n+1-i,-y)$, where $R'$ might be empty but where $(1,n+1-i,-y)$ contains at least two vertices. If Right can take fully $G_1$, he takes it and at the end, the score will be at most 4.
Otherwise, we must have $n\geq 2$ and then we have $i\leq n$. Right plays in column $i$ of $G_1$ and take as many vertices as Left. Then it is remaining after the two moves $R+(1,n-i-1,y)+R'+(-1,n-i+1,-y)$ and we can conclude as before.

Finally, assume that Left takes a full component. If Right can take the other component, the game ends with a score of at least 4.
Otherwise, it means that Left played in $G_1$ and that $n\leq 3$. Let $i\leq 2$ be the column where Left played.
If $n=1$, the only case where Right cannot take the whole $G_2$ is when $G_2$ has form $(1,3,1)$ or $(0,3,1)$. In the first case, Right can take four vertices and then both players take two vertices, ending in a draw. In the second case, Right takes five vertices, Left can take the two remaining vertices, and the games end in a draw.
If $n=2$, Left took at most six vertices. If $x$ or $y$ is equal to $-1$, then Right can take four vertices in $G_2$ and leaves a graph with two components of size at least 2 by playing in column $2$ or $3$. Thus at the end, Right will win at least six vertices. If $x=y=1$, Right can directly take six vertices in $G_2$ by playing in column 2. In both cases, Left takes at most ten vertices, ending with a score of at most 4. If $x=0$ or $y=0$, Right can take 5 vertices and there are at most 14 vertices in total, thus the score is also at most 4. 
Finally, if $n=3$, Left took at most eight vertices. 
As before, Right can take at least $|G_1|$ vertices by taking four vertices in his first move and splitting $G_2$ in two components of size at most 3 (if $|G_1|=7$) or $4$ (if $|G_1|=8$).



Turn back to our original graph $G_{2,m}$. Right can follow the following strategy :
\begin{itemize}
    \item While Left does not play in the five columns in the middle, play the symmetric. The most left and most right graphs are opposite and thus simplify each other. The remaining graph has form $(1,m',-1)$ with $m'$ odd. In particular, there is still a symmetry center.
    \item If Left plays column $(m+1)/2-2$ (or the symmetric), play the symmetric. The components on the left and the right simplify. Only four vertices are left in the middle, so Left win at most four vertices. We have $Ls(G_{2,m}) \le 4$.
    \item If Left plays at distance 1 to the central column, i.e. she plays $(m+1)/2-1$, then play the symmetric. The two remaining graphs are opposite. So $Ls(G_{2,m}) \le 0$.
    \item If Left plays the central column $(m+1)/2$:
    \begin{itemize}
        \item If the vertices in column $(m+1)/2-3$ have been played, the graph was (before the move of Left) of the form $(1,3,-1)$. Thus Left took six vertices and let a $S_2$. Right can take it. The final score is $6-2 = 4$.
        \item If the vertices in column $(m+1)/2-3$ have not been played, Right plays in column $(m+1)/2-2$ and takes four vertices (like Left did). Then the remaining graph has form $(x,n,1) + (-1,n+2,-x)$. Using the result proved above (recall that $(-1,n+2,-x)=(-x,n+2,-1)$, the score is at most 4 at the end.
    \end{itemize}
\end{itemize}

Finally, $Ls(G_{2,m}) \le 4$, and the theorem is proven. 
\end{proof}

\begin{proof}[Proof of Lemma \ref{lem:2rowseven}]
Let $m\geq 4$ be an even integer.  We prove that $Ls(G_{2,m})\geq 2$ by giving a strategy for Left to score $2$. Left can start by playing in column $m/2$ (it corresponds to the vertex marked by a circle on Figure \ref{fig:even}). Denote by $G'$ the remaining graph after this first move. Left can then decide to "give" the two vertices marked by blue triangles in the graph of Figure \ref{fig:even}. More precisely, we apply Theorem \ref{thmineq} with $B_0 = \{  x_1, x_2\}$ where $x_1$ and $x_2$ are the black vertices of columns $m/2+1$ and $m/2+2$. We obtain $ Rs(G') \ge  Rs(G'\bs B_0) - |B_0| \ge -2$, since $G'\bs B_0$ is the sum of a game and its opposite, and thus a score $0$. Finally, $Ls(G_{2,m}) \ge 4- Rs(G') \ge 2$.

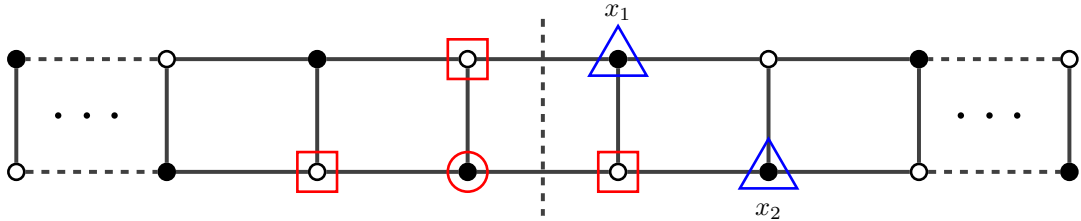
\begin{figure}[ht]
    \centering

\begin{tikzpicture}

\Vertex[x=1,y=-.75, style = {color = white}, label = \dots, size = .5, fontscale = 3]{16}

\Vertex[x=13, y=-.75, style = {color = white}, label = \dots, size = .5, fontscale = 3]{16}

    \draw (0,0) node[vertexB](1){} ;

    \draw (0,-1.5) node[vertexW](2){} ;

    \draw (2,-1.5) node[vertexB](3){} ;

    \draw (2,0) node[vertexW](4){} ;

    \draw (4,0) node[vertexB](5){} ;
  \draw (4,-1.5) node[vertexW](6){} ;

  \draw (4,-1.5) node[redsquare](n){} ;
  
  \draw (6,-1.5) node[vertexB](7){} ;
    
  \draw (6,-1.5) node[redcircle](n){} ;

  \draw (6,0) node[vertexW](8){} ;
  
  \draw (6,0) node[redsquare](n){} ;
  
  \draw (8,0) node[vertexB](9){} ;

  \draw (8,0) node[bluetriangle](n1){} node[above=0.4] {$x_1$} ;
  
  \draw (8,-1.5) node[vertexW](10){} ;

  \draw (8,-1.5) node[redsquare](n){} ;

  \draw (10,-1.5) node[vertexB](11){} ;
  \draw (10,-1.5) node[bluetriangle](n){} node[below=0.3] {$x_2$} ;

  \draw (10,0) node[vertexW](12){} ;

  \draw (12,0) node[vertexB](13){} ;

  \draw (12,-1.5) node[vertexW](14){} ;

  \draw (14,-1.5) node[vertexB](15){} ;

  \draw (14,0) node[vertexW](16){} ;

  \Vertex[x = 7, y = -2.2, style = {color = white}, size = .001]{20}
  \Vertex[x = 7, y = .7, style = {color = white}, size = .001]{21}

  \Edge[, style = dashed](4)(1)
  \Edge[](2)(1)
  \Edge[, style = dashed](2)(3)
  \Edge[](4)(3) 
  \Edge[](4)(5)
  \Edge[](6)(3)
  \Edge[](6)(5)
  \Edge[](6)(7)
  \Edge[](8)(5)
  \Edge[](8)(7)
  \Edge[](8)(9)
  \Edge[](10)(7)
  \Edge[](10)(9)
  \Edge[](10)(11)
  \Edge[](12)(9)
  \Edge[](12)(11)
  \Edge[](12)(13)
  \Edge[](14)(11)
  \Edge[](14)(13)
  \Edge[, style = dashed](14)(15)
  \Edge[, style = dashed](16)(13)
  \Edge[](16)(15)
  \Edge[style = dashed](20)(21)

\end{tikzpicture}

    \caption{Strategy if the number of columns is even}
    \label{fig:even}
\end{figure}

\end{proof}

Theorem \ref{thm:2rowsodd} and Lemma \ref{lem:2rowseven} actually suggest the following conjecture that gives the complete values for grids of two rows. It only remains to prove that $Ls(G_{2,m})\leq 2$ for even $m$.
We have checked this conjecture up to $m=30$.

\begin{conjecture}
Let $m\geq 2$. We have 

$$Ls(G_{2,m}) = -Rs(G_{2,m}) = \left\{
    \begin{array}{lll}
        6 & \mbox{if } m = 3;  \\
        4 & \mbox{if } m \equiv 1 \mod 2 \mbox{ and } m \ge 5 \mbox{ or } m = 2;  \\
        2 & \mbox{if } m \equiv 0 \mod 2 \mbox{ and } m \ge 4.
    \end{array}
\right.$$

\end{conjecture}

\subsection{Grids of three rows}

As detailed in the previous sections, in the cases of segments and grids of two rows, a common technique to prove that the first player has a winning strategy consists in cutting the graph into two "similar" parts and then guarantee that he will not lose many points on the rest of the game. Unfortunately, this technique is no more available for grids having at least four rows, as single moves generally keep the graph connected. Yet, in the case of grids of three rows, some moves split the grid into two parts, but other moves keep it connected. This is sufficient to yield strategies for the first player in almost all configurations.\\

First note that if $m$ is odd, as the grid is not symmetric, having $Ls(G_{3,m}) \ge 0$ or $Rs(G_{3,m}) \le 0$ is not straightforward at first sight. However, we will prove Conjecture~\ref{conj:grids} for grids of three rows in almost all cases.

\begin{theorem}
Let $G_{3,m}$ be a grid of three rows and $m$ columns, with a black vertex at the top left. We have $Ls(G_{3,m}) > 0$, and if $m \not \equiv 3 \mod 4$,  $Rs(G_{3,m}) < 0$.
\end{theorem}

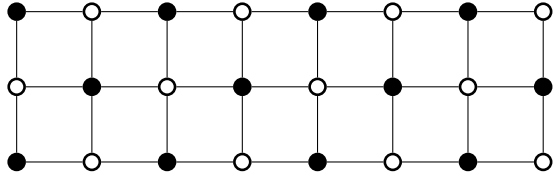
\begin{figure}
    \centering
    
\begin{tikzpicture}

  \draw (0,0) node[vertexW](1){} ;
  \draw (0,1) node[vertexB](2){} ;
  \draw (1,0) node[vertexB](3){} ;
  \draw (1,1) node[vertexW](4){} ;
  \draw (2,0) node[vertexW](5){} ;
  \draw (2,1) node[vertexB](6){} ;
  \draw (3,0) node[vertexB](7){} ;
  \draw (3,1) node[vertexW](8){} ;
  \draw (4,0) node[vertexW](9){} ;
  \draw (4,1) node[vertexB](10){} ;
  \draw (5,0) node[vertexB](11){} ;
  \draw (5,1) node[vertexW](12){} ;
  \draw (6,0) node[vertexW](13){} ;
  \draw (6,1) node[vertexB](14){} ;
  \draw (7,0) node[vertexB](15){} ;
  \draw (7,1) node[vertexW](16){} ;

  \draw (0,-1) node[vertexB](17){} ;
  \draw (1,-1) node[vertexW](18){} ;
  \draw (2,-1) node[vertexB](19){} ;
  \draw (3,-1) node[vertexW](20){} ;
  \draw (4,-1) node[vertexB](21){} ;
  \draw (5,-1) node[vertexW](22){} ;
  \draw (6,-1) node[vertexB](23){} ;
  \draw (7,-1) node[vertexW](24){} ;

  \path[draw,-] (1) to (17);
  \path[draw,-] (3) to (18);
  \path[draw,-] (5) to (19);
  \path[draw,-] (7) to (20);
  \path[draw,-] (9) to (21);
  \path[draw,-] (11) to (22);
  \path[draw,-] (13) to (23);
  \path[draw,-] (15) to (24);

  \path[draw,-] (17) to (18);
  \path[draw,-] (18) to (19);
  \path[draw,-] (19) to (20);
  \path[draw,-] (20) to (21);
  \path[draw,-] (21) to (22);
  \path[draw,-] (22) to (23);
  \path[draw,-] (23) to (24);
  
  \path[draw,-] (1) to (2);
  \path[draw,-] (1) to (3);
  \path[draw,-] (2) to (4);
  \path[draw,-] (3) to (4);
  \path[draw,-] (3) to (5);
  \path[draw,-] (4) to (6);
  \path[draw,-] (5) to (6);
  \path[draw,-] (5) to (7);
  \path[draw,-] (6) to (8);
  \path[draw,-] (7) to (8);
  \path[draw,-] (7) to (9);
  \path[draw,-] (8) to (10);
  \path[draw,-] (9) to (10);
  \path[draw,-] (9) to (11);
  \path[draw,-] (10) to (12);
  \path[draw,-] (11) to (12);
  \path[draw,-] (11) to (13);
  \path[draw,-] (12) to (14);
  \path[draw,-] (13) to (14);
  \path[draw,-] (13) to (15);
  \path[draw,-] (14) to (16);
  \path[draw,-] (15) to (16);

\end{tikzpicture} 
    
    \caption{Example of grid of three rows and eight columns.}
    \label{fig:G3}
\end{figure}

\begin{proof}

To make the proof more understandable, in the following figures, we will mark  by a red circle the vertex played by the first player, by red squares the vertices taken by this move, and by blue triangles the vertices freely given to the opponent using Theorem \ref{thmineq}.\\

Let $ G = G_{3,m}$ be a grid of three rows. 

\begin{itemize}
    \item If $m \equiv 0 \mod 2$. We already know by symmetry that $Ls(G) = - Rs(G)$. So, without loss of generality, it is sufficient to prove that $Ls(G) > 0$.
    
    Left can adopt the following strategy: play the black vertex on column $\frac{m}{2}$ marked by a red circle (see Figure~\ref{fig:G30}). Left takes all vertices marked by red squares and scores $5$ points. Then, she can give the two black vertices $x_1$ and $x_2$ of column $\frac{m}{2} + 1$ to Right, and apply Theorem \ref{thmineq} with $B_0 = \{x_1, x_2\}$.
    
    The remaining graph is composed of two grids $G_{3,\frac{m}{2} -1}$. But the left part has one missing white vertex marked by a red square square. Add it to the graph by giving one more point to Right once again according to Theorem \ref{thmineq} with $W_0 = \{x_3\}$. Note $H_{sym}$ the symmetric graph obtained. The final score thus satisfies $Ls(G) \ge 5 + Rs(H_{sym}) - |W_0| - |B_0| \ge 5 - 1 - 2 = 2$  .

\begin{figure}[ht]
    \centering
\begin{tikzpicture}

  \draw (0,0) node[vertexW](1){} ;
  \draw (0,1) node[vertexB](2){} ;
  \draw (1,0) node[vertexB](3){} ;
  \draw (1,1) node[vertexW](4){} ;
  \draw (2,0) node[vertexW](5){} ;
  \draw (2,1) node[vertexB](6){} ;
  \draw (3,0) node[vertexB](7){} ;
  \draw (3,1) node[vertexW](8){} ;
  \draw (4,0) node[vertexW](9){} ;
  \draw (4,1) node[vertexB](10){} ;
  \draw (5,0) node[vertexB](11){} ;
  \draw (5,1) node[vertexW](12){} ;
  \draw (6,0) node[vertexW](13){} ;
  \draw (6,1) node[vertexB](14){} ;
  \draw (7,0) node[vertexB](15){} ;
  \draw (7,1) node[vertexW](16){} ;

  \draw (0,-1) node[vertexB](17){} ;
  \draw (1,-1) node[vertexW](18){} ;
  \draw (2,-1) node[vertexB](19){} ;
  \draw (3,-1) node[vertexW](20){} ;
  \draw (4,-1) node[vertexB](21){} ;
  \draw (5,-1) node[vertexW](22){} ;
  \draw (6,-1) node[vertexB](23){} ;
  \draw (7,-1) node[vertexW](24){} ;

  \path[draw,-] (1) to (17);
  \path[draw,-] (3) to (18);
  \path[draw,-] (5) to (19);
  \path[draw,-] (7) to (20);
  \path[draw,-] (9) to (21);
  \path[draw,-] (11) to (22);
  \path[draw,-] (13) to (23);
  \path[draw,-] (15) to (24);

  \path[draw,-] (17) to (18);
  \path[draw,-] (18) to (19);
  \path[draw,-] (19) to (20);
  \path[draw,-] (20) to (21);
  \path[draw,-] (21) to (22);
  \path[draw,-] (22) to (23);
  \path[draw,-] (23) to (24);
  
  \path[draw,-] (1) to (2);
  \path[draw,-] (1) to (3);
  \path[draw,-] (2) to (4);
  \path[draw,-] (3) to (4);
  \path[draw,-] (3) to (5);
  \path[draw,-] (4) to (6);
  \path[draw,-] (5) to (6);
  \path[draw,-] (5) to (7);
  \path[draw,-] (6) to (8);
  \path[draw,-] (7) to (8);
  \path[draw,-] (7) to (9);
  \path[draw,-] (8) to (10);
  \path[draw,-] (9) to (10);
  \path[draw,-] (9) to (11);
  \path[draw,-] (10) to (12);
  \path[draw,-] (11) to (12);
  \path[draw,-] (11) to (13);
  \path[draw,-] (12) to (14);
  \path[draw,-] (13) to (14);
  \path[draw,-] (13) to (15);
  \path[draw,-] (14) to (16);
  \path[draw,-] (15) to (16);

  \draw (3,0) node[redcircle](n1){} ;
  \draw (4,0) node[redsquare](n1){} ;
  \draw (2,0) node[redsquare](n1){} node[below left=0.2] {$x_3$} ;
  \draw (3,-1) node[redsquare](n1){} ;
  \draw (3,1) node[redsquare](n1){} ;
  \draw (4,1) node[bluetriangle](n1){} node[above=0.4] {$x_1$} ;
  \draw (4,-1) node[bluetriangle](n1){} node[below=0.3] {$x_2$};

\end{tikzpicture}     \caption{Strategy if $m \equiv 0 \mod 2$}
    \label{fig:G30}
\end{figure}
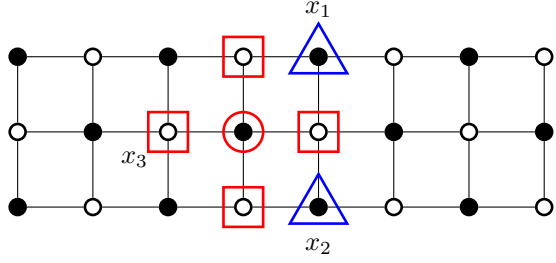  

    \item If $m \equiv 1 \mod 4$ and Left starts, she can play the black vertex at the bottom of column $\frac{m+1}{2}$ marked by a red circle (see Figure  \ref{fig:G31L}). Then, by giving to Right the two black vertices in top corners and the black vertex in the middle of the first row, she creates two identical symmetric graphs. Note this set of three vertices $B_0 = \{x_1, x_2, x_3 \}$ and the resulting graph $H_{sym}$. Therefore, the score satisfies $Ls(G) \ge 4 + Rs(H_{sym}) - |B_0| = 4 + 0 -3 = 1$.

\begin{figure}[ht]
    \centering
\begin{tikzpicture}

  \draw (0,0) node[vertexW](1){} ;
  \draw (0,1) node[vertexB](2){} ;
  \draw (1,0) node[vertexB](3){} ;
  \draw (1,1) node[vertexW](4){} ;
  \draw (2,0) node[vertexW](5){} ;
  \draw (2,1) node[vertexB](6){} ;
  \draw (3,0) node[vertexB](7){} ;
  \draw (3,1) node[vertexW](8){} ;
  \draw (4,0) node[vertexW](9){} ;
  \draw (4,1) node[vertexB](10){} ;
  \draw (5,0) node[vertexB](11){} ;
  \draw (5,1) node[vertexW](12){} ;
  \draw (6,0) node[vertexW](13){} ;
  \draw (6,1) node[vertexB](14){} ;
  \draw (7,0) node[vertexB](15){} ;
  \draw (7,1) node[vertexW](16){} ;
  \draw (8,0) node[vertexW](25){} ;
  \draw (8,1) node[vertexB](26){} ;

  \draw (0,-1) node[vertexB](17){} ;
  \draw (1,-1) node[vertexW](18){} ;
  \draw (2,-1) node[vertexB](19){} ;
  \draw (3,-1) node[vertexW](20){} ;
  \draw (4,-1) node[vertexB](21){} ;
  \draw (5,-1) node[vertexW](22){} ;
  \draw (6,-1) node[vertexB](23){} ;
  \draw (7,-1) node[vertexW](24){} ;
  \draw (8,-1) node[vertexB](27){} ;

  \path[draw,-] (1) to (17);
  \path[draw,-] (3) to (18);
  \path[draw,-] (5) to (19);
  \path[draw,-] (7) to (20);
  \path[draw,-] (9) to (21);
  \path[draw,-] (11) to (22);
  \path[draw,-] (13) to (23);
  \path[draw,-] (15) to (24);
  \path[draw,-] (16) to (26);
  \path[draw,-] (15) to (25);
  \path[draw,-] (24) to (27);
  \path[draw,-] (27) to (25);
  \path[draw,-] (25) to (26);

  \path[draw,-] (17) to (18);
  \path[draw,-] (18) to (19);
  \path[draw,-] (19) to (20);
  \path[draw,-] (20) to (21);
  \path[draw,-] (21) to (22);
  \path[draw,-] (22) to (23);
  \path[draw,-] (23) to (24);
  
  \path[draw,-] (1) to (2);
  \path[draw,-] (1) to (3);
  \path[draw,-] (2) to (4);
  \path[draw,-] (3) to (4);
  \path[draw,-] (3) to (5);
  \path[draw,-] (4) to (6);
  \path[draw,-] (5) to (6);
  \path[draw,-] (5) to (7);
  \path[draw,-] (6) to (8);
  \path[draw,-] (7) to (8);
  \path[draw,-] (7) to (9);
  \path[draw,-] (8) to (10);
  \path[draw,-] (9) to (10);
  \path[draw,-] (9) to (11);
  \path[draw,-] (10) to (12);
  \path[draw,-] (11) to (12);
  \path[draw,-] (11) to (13);
  \path[draw,-] (12) to (14);
  \path[draw,-] (13) to (14);
  \path[draw,-] (13) to (15);
  \path[draw,-] (14) to (16);
  \path[draw,-] (15) to (16);

  \draw (4,-1) node[redcircle](n1){} ;
  \draw (4,0) node[redsquare](n1){} ;
  \draw (3,-1) node[redsquare](n1){} ;
  \draw (5,-1) node[redsquare](n1){} ;
  \draw (4,1) node[bluetriangle](n1){} node[above =0.4] {$x_2$};
  \draw (0,1) node[bluetriangle](n1){} node[above=0.4] {$x_1$};
  \draw (8,1) node[bluetriangle](n1){} node[above =0.4] {$x_3$};

\end{tikzpicture}     \caption{Strategy for Left if $m \equiv 1 \mod 4$.}
    \label{fig:G31L}
\end{figure}
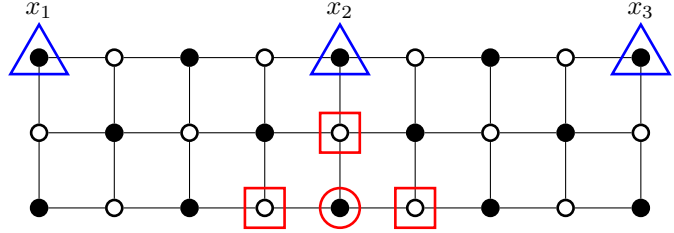

    \item If $m \equiv 1 \mod 4$ and Right starts, he can play the vertex in the middle of the graph (see Figure \ref{fig:G31R}), and then give the two black vertices $x_1$ and $x_2$ of the columns $\frac{m-1}{2}$ and $\frac{m+3}{2}$ to Left. Note them $B_0 = \{x_1, x_2\}$. The remaining graph is composed of two instances of $G_{3,\frac{m-1}{2}}$, and is symmetric as $m-1 \equiv 0 \mod 4$. Note it $H_{sym} = 2*G_{3,\frac{m-1}{2}}$. Thus, by Theorem \ref{thmineq}, the score verifies $Rs(G) \le -5 + Ls(H_{sym}) + |B_0| =-5 + 0 +2 = -3$

\begin{figure}[ht]
    \centering
\begin{tikzpicture}

  \draw (0,0) node[vertexW](1){} ;
  \draw (0,1) node[vertexB](2){} ;
  \draw (1,0) node[vertexB](3){} ;
  \draw (1,1) node[vertexW](4){} ;
  \draw (2,0) node[vertexW](5){} ;
  \draw (2,1) node[vertexB](6){} ;
  \draw (3,0) node[vertexB](7){} ;
  \draw (3,1) node[vertexW](8){} ;
  \draw (4,0) node[vertexW](9){} ;
  \draw (4,1) node[vertexB](10){} ;
  \draw (5,0) node[vertexB](11){} ;
  \draw (5,1) node[vertexW](12){} ;
  \draw (6,0) node[vertexW](13){} ;
  \draw (6,1) node[vertexB](14){} ;
  \draw (7,0) node[vertexB](15){} ;
  \draw (7,1) node[vertexW](16){} ;
  \draw (8,0) node[vertexW](25){} ;
  \draw (8,1) node[vertexB](26){} ;

  \draw (0,-1) node[vertexB](17){} ;
  \draw (1,-1) node[vertexW](18){} ;
  \draw (2,-1) node[vertexB](19){} ;
  \draw (3,-1) node[vertexW](20){} ;
  \draw (4,-1) node[vertexB](21){} ;
  \draw (5,-1) node[vertexW](22){} ;
  \draw (6,-1) node[vertexB](23){} ;
  \draw (7,-1) node[vertexW](24){} ;
  \draw (8,-1) node[vertexB](27){} ;

  \path[draw,-] (1) to (17);
  \path[draw,-] (3) to (18);
  \path[draw,-] (5) to (19);
  \path[draw,-] (7) to (20);
  \path[draw,-] (9) to (21);
  \path[draw,-] (11) to (22);
  \path[draw,-] (13) to (23);
  \path[draw,-] (15) to (24);
  \path[draw,-] (16) to (26);
  \path[draw,-] (15) to (25);
  \path[draw,-] (24) to (27);
  \path[draw,-] (27) to (25);
  \path[draw,-] (25) to (26);

  \path[draw,-] (17) to (18);
  \path[draw,-] (18) to (19);
  \path[draw,-] (19) to (20);
  \path[draw,-] (20) to (21);
  \path[draw,-] (21) to (22);
  \path[draw,-] (22) to (23);
  \path[draw,-] (23) to (24);
  
  \path[draw,-] (1) to (2);
  \path[draw,-] (1) to (3);
  \path[draw,-] (2) to (4);
  \path[draw,-] (3) to (4);
  \path[draw,-] (3) to (5);
  \path[draw,-] (4) to (6);
  \path[draw,-] (5) to (6);
  \path[draw,-] (5) to (7);
  \path[draw,-] (6) to (8);
  \path[draw,-] (7) to (8);
  \path[draw,-] (7) to (9);
  \path[draw,-] (8) to (10);
  \path[draw,-] (9) to (10);
  \path[draw,-] (9) to (11);
  \path[draw,-] (10) to (12);
  \path[draw,-] (11) to (12);
  \path[draw,-] (11) to (13);
  \path[draw,-] (12) to (14);
  \path[draw,-] (13) to (14);
  \path[draw,-] (13) to (15);
  \path[draw,-] (14) to (16);
  \path[draw,-] (15) to (16);

  \draw (4,0) node[redcircle](n1){} ;
  \draw (4,1) node[redsquare](n1){} ;
  \draw (4,-1) node[redsquare](n1){} ;
  \draw (3,0) node[redsquare](n1){} node[below left=0.2] {$x_1$};
  \draw (5,0) node[redsquare](n1){} node[below right=0.2] {$x_2$};

\end{tikzpicture}     \caption{Strategy for Right if $m \equiv 1 \mod 4$.}
    \label{fig:G31R}
\end{figure}
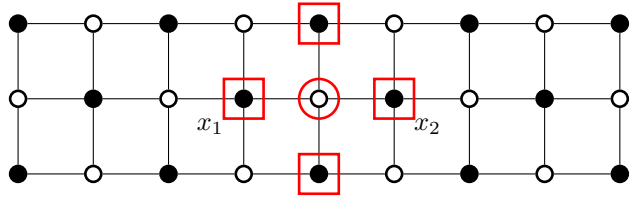 
    
\item If $m \equiv 3 \mod 4$ and Left starts, she can play the middle vertex of the graph and give the four black vertices at its diagonals to Right (see Figure \ref{fig:G33L}). Note these four vertices $B_0 = \{x_1, x_2, x_3, x_4\}$. The remaining graph $H_{sym}$ is composed of two instances of $G_{3,\frac{m-3}{2}}$, so it is equivalent to $0$ as $m-3 \equiv 0 \mod 4$. Finally, the score satisfies $Ls(G) \ge 5 + Rs(H_{sym}) - |B_0| = 5 + 0 - 4 = 1$.

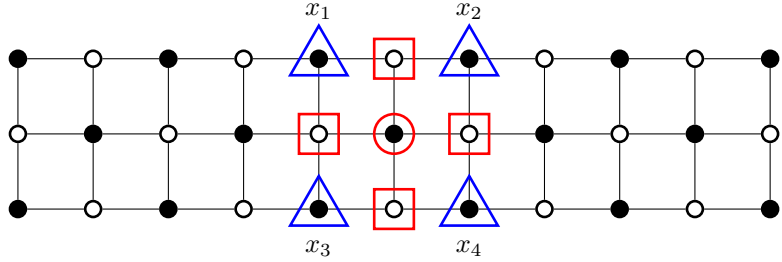
\begin{figure}[ht]
    \centering
\begin{tikzpicture}

  \draw (0,0) node[vertexW](1){} ;
  \draw (0,1) node[vertexB](2){} ;
  \draw (1,0) node[vertexB](3){} ;
  \draw (1,1) node[vertexW](4){} ;
  \draw (2,0) node[vertexW](5){} ;
  \draw (2,1) node[vertexB](6){} ;
  \draw (3,0) node[vertexB](7){} ;
  \draw (3,1) node[vertexW](8){} ;
  \draw (4,0) node[vertexW](9){} ;
  \draw (4,1) node[vertexB](10){} ;
  \draw (5,0) node[vertexB](11){} ;
  \draw (5,1) node[vertexW](12){} ;
  \draw (6,0) node[vertexW](13){} ;
  \draw (6,1) node[vertexB](14){} ;
  \draw (7,0) node[vertexB](15){} ;
  \draw (7,1) node[vertexW](16){} ;
  \draw (8,0) node[vertexW](25){} ;
  \draw (8,1) node[vertexB](26){} ;
  \draw (9,0) node[vertexB](27){} ;
  \draw (9,1) node[vertexW](28){} ;
  \draw (10,0) node[vertexW](29){} ;
  \draw (10,1) node[vertexB](30){} ;

  \draw (0,-1) node[vertexB](17){} ;
  \draw (1,-1) node[vertexW](18){} ;
  \draw (2,-1) node[vertexB](19){} ;
  \draw (3,-1) node[vertexW](20){} ;
  \draw (4,-1) node[vertexB](21){} ;
  \draw (5,-1) node[vertexW](22){} ;
  \draw (6,-1) node[vertexB](23){} ;
  \draw (7,-1) node[vertexW](24){} ;
  \draw (8,-1) node[vertexB](31){} ;
  \draw (9,-1) node[vertexW](32){} ;
  \draw (10,-1) node[vertexB](33){} ;

  \path[draw,-] (1) to (17);
  \path[draw,-] (3) to (18);
  \path[draw,-] (5) to (19);
  \path[draw,-] (7) to (20);
  \path[draw,-] (9) to (21);
  \path[draw,-] (11) to (22);
  \path[draw,-] (13) to (23);
  \path[draw,-] (15) to (24);
  \path[draw,-] (16) to (26);
  \path[draw,-] (15) to (25);
  \path[draw,-] (24) to (31);
  \path[draw,-] (31) to (25);
  \path[draw,-] (25) to (26);
  
  \path[draw,-] (27) to (28);
  \path[draw,-] (27) to (32);
  \path[draw,-] (29) to (30);
  \path[draw,-] (29) to (33);
  \path[draw,-] (25) to (27);
  \path[draw,-] (27) to (29);

  \path[draw,-] (26) to (28);
  \path[draw,-] (28) to (30);
  \path[draw,-] (31) to (32);
  \path[draw,-] (32) to (33);

  \path[draw,-] (17) to (18);
  \path[draw,-] (18) to (19);
  \path[draw,-] (19) to (20);
  \path[draw,-] (20) to (21);
  \path[draw,-] (21) to (22);
  \path[draw,-] (22) to (23);
  \path[draw,-] (23) to (24);
  
  \path[draw,-] (1) to (2);
  \path[draw,-] (1) to (3);
  \path[draw,-] (2) to (4);
  \path[draw,-] (3) to (4);
  \path[draw,-] (3) to (5);
  \path[draw,-] (4) to (6);
  \path[draw,-] (5) to (6);
  \path[draw,-] (5) to (7);
  \path[draw,-] (6) to (8);
  \path[draw,-] (7) to (8);
  \path[draw,-] (7) to (9);
  \path[draw,-] (8) to (10);
  \path[draw,-] (9) to (10);
  \path[draw,-] (9) to (11);
  \path[draw,-] (10) to (12);
  \path[draw,-] (11) to (12);
  \path[draw,-] (11) to (13);
  \path[draw,-] (12) to (14);
  \path[draw,-] (13) to (14);
  \path[draw,-] (13) to (15);
  \path[draw,-] (14) to (16);
  \path[draw,-] (15) to (16);

  \draw (5,0) node[redcircle](n1){} ;
  \draw (5,1) node[redsquare](n1){} ;
  \draw (5,-1) node[redsquare](n1){} ;
  \draw (4,0) node[redsquare](n1){} ;
  \draw (6,0) node[redsquare](n1){} ;
  \draw (4,1) node[bluetriangle](n1){} node[above =0.4]{$x_1$} ;
  \draw (4,-1) node[bluetriangle](n1){} node[below =0.3]{$x_3$};
  \draw (6,1) node[bluetriangle](n1){} node[above =0.4]{$x_2$};
  \draw (6,-1) node[bluetriangle](n1){} node[below =0.3]{$x_4$};

\end{tikzpicture}     \caption{Strategy for Left if $m \equiv 3 \mod 4$.}
    \label{fig:G33L}
\end{figure}

\end{itemize}

\end{proof}

The case of grids of three rows is almost solved from the above result. As an open problem, it remains to show that if $m \equiv 3 \mod 4$ and Right starts, he has a winning strategy. 

\bibliography{Bipartite_instances_of_Influence}   
\bibliographystyle{alpha} 

\end{document}